\documentclass[reqno,11pt]{amsart}
\usepackage{graphicx, amsmath, amssymb, amscd, mathtools, hyperref, amsthm, euscript, 
amsfonts,bm,color}
\hypersetup{hypertexnames=false}
\usepackage{graphicx}
 
\usepackage{mathtools}
\usepackage{amsmath}
\usepackage{amssymb}
 
\usepackage{amsthm}
\usepackage{graphicx}
\usepackage{mathrsfs}
\usepackage{xcolor}
\usepackage{tikz-cd}
\usepackage{tikz}
\usetikzlibrary{patterns}
\usepackage{hyperref}
\usepackage{cleveref}
 
\usetikzlibrary{arrows.meta}

\usepackage{verbatim}
\usepackage{float}
 
\usepackage{enumitem}

\usepackage{thmtools}
\usepackage{thm-restate}
 
\usepackage{caption}

\DeclareMathAlphabet{\mathpzc}{OT1}{pzc}{m}{it}
 
\numberwithin{equation}{section}

\newtheorem{theorem}{Theorem}[section]
 
\newtheorem{corollary}[theorem]{Corollary}

\newtheorem{proposition}[theorem]{Proposition}

\newtheorem{Remark}[theorem]{Remark}
 
\newtheorem{lemma}[theorem]{Lemma}
 
\theoremstyle{definition}

\newtheorem{Def}[theorem]{Definition}

\newcommand{\be}{begin{equation}}

\newcommand{\e}{\epsilon}

\newcommand{{\grinv}}{{\Cal G}^{-r}}

\newcommand{\Cal}{\mathcal}

\newcommand{\bp}{\begin{pmatrix}}
\newcommand{\ep}{\end{pmatrix}}
\renewcommand{\be}{\begin{equation}}
\newcommand{\ee}{\end{equation}}
\renewcommand{\epsilon}{\varepsilon}

\renewcommand{\bp}{{\rm bp}}

\newcommand{\T}{\operatorname{T}}

\newcommand{\norm}[1]{\lVert #1 \rVert}

\newcommand{\supp}{\operatorname{supp}}

\newcommand{\BMS}{\operatorname{BMS}}

\renewcommand{\setminus}{-}

\newcommand{\ga}{\gamma}
\newcommand{\F}{\mathcal F}
\newcommand{\La}{\Lambda}
\renewcommand{\i}{\op{i}}

\def\scrL{{\mathcal L}}

\def\i{\mathrm{i}}

\def\dim{\operatorname{dim}}

\def\T{\operatorname{T{}}}

\def\supp{\operatorname{supp}}

\newcommand{\sfrac}[2]{{\textstyle \frac {#1}{#2}}}

\newcommand{\fa}{\mathfrak{a}}

\newcommand{\Ga}{\Gamma}
\newcommand{\cal}{\mathcal}
\renewcommand{\e}{\varepsilon}

\def\scrL{{\mathcal L}}

\def\i{\mathrm{i}}

\def\dim{\operatorname{dim}}

\def\T{\operatorname{T{}}}

\def\supp{\operatorname{supp}}

\def\supp{\operatorname{supp}}

\newcommand{\te}{\theta}

\newcommand{\err}{\mathsf{E}}

\newcommand{\gap}{{\mathsf m}_\theta}
 
\begin{document}
 
\title[Local limit theorems and counting]{Local limit theorems and counting for transverse subgroups}

\begin{abstract} 
Transverse subgroups form a special class of discrete subgroups of connected semisimple real algebraic groups that includes all discrete subgroups of rank-one Lie groups, as well as Anosov and relatively Anosov subgroups.
We study the joint distribution of partial Cartan projections and the two associated Cartan flags for Zariski dense $\theta$-transverse subgroups.
Assuming local mixing of the associated Bowen-Margulis-Sullivan measure, we prove two local limit theorems and deduce counting asymptotics in Euclidean balls and tubes, together with equidistribution of the Cartan flags. 
These qualitative conclusions hold for $\theta$-Anosov subgroups, for which the local mixing hypothesis is satisfied.
If the group is linearly $\theta$-regular near the direction of maximal growth, a polynomial error in local mixing gives polynomial errors in these results.
Applying the quantitative local mixing theorem of Kim-Oh-Pan, we obtain local limit and counting theorems with polynomial error estimates for several classes of relatively Anosov representations of Fuchsian groups.
\end{abstract}

\author{Sam Edwards}
\address{Department of Mathematics, Uppsala University, Box 524, 751 20 Uppsala, Sweden}
\email{samuel.edwards@math.uu.se}
 
\author{Minju Lee}
\address{Department of Mathematical Sciences, KAIST, 291 Daehak-ro, Yuseong-gu, Daejeon 34141, Republic of Korea}
\email{minju.lee@kaist.ac.kr}
 
\author{Hee Oh}
\address{Department of Mathematics, Yale University, PO Box 208283, New Haven, CT 06520-8283, USA}
\email{hee.oh@yale.edu}
\thanks{Oh is partially supported by NSF grant DMS-2450703.}
 
\maketitle
\tableofcontents

\section{Introduction}
Transverse subgroups form a special class of discrete subgroups of connected semisimple real algebraic groups that includes all discrete subgroups of rank-one Lie groups, as well as Anosov and relatively Anosov subgroups.
For Zariski dense $\theta$-transverse subgroups, we study the joint distribution of partial Cartan projections and the two associated Cartan flags.
Assuming local mixing of the associated Bowen-Margulis-Sullivan measure in the direction of maximal growth, with transverse perturbations, we prove two local limit theorems and deduce counting asymptotics in Euclidean balls and tubes. Under a linear regularity assumption, a polynomial error in local mixing gives polynomial errors in these results.

 \medskip 
To state our results, let $G$ be a connected semisimple real algebraic group. Fix a Cartan decomposition $G=K\exp(\mathfrak a^+)K$, let $\Pi$ be the corresponding set of simple roots, and fix a nonempty subset $\theta\subset\Pi$. 
We assume that $$\text{$\Gamma<G$ is a Zariski dense $\theta$-transverse subgroup}$$ (see Section \ref{ps}). This class includes Zariski dense subgroups of relatively $\theta$-Anosov subgroups. Relatively $\theta$-Anosov subgroups generalize geometrically finite subgroups of rank-one Lie groups and include $\theta$-Anosov subgroups. Anosov subgroups were introduced by Labourie \cite{Lab} and extended to word hyperbolic groups by Guichard-Wienhard \cite{GW}. For the relative theory, see Kapovich-Leeb \cite{KL} and Zhu-Zimmer \cite{ZZ1}.

The subset $\theta$ determines the Cartan subspace and flag varieties used in our counting problem. 
Write
$$
 \mathfrak a_\theta=\bigcap_{\alpha\in\Pi\setminus\theta}\ker\alpha,\qquad \mathfrak a_\theta^+=\mathfrak a_\theta\cap\mathfrak a^+,\qquad A_\theta=\exp\mathfrak a_\theta,
$$
and let $A_\theta^+=\exp\mathfrak a_\theta^+$. Fix a Weyl-group invariant inner product on $\mathfrak a$, with associated norm $\|\cdot\|$, and let $p_\theta:\mathfrak a\to\mathfrak a_\theta$ be the orthogonal projection.

For $g\in G$, its Cartan projection $\mu(g)\in\mathfrak a^+$ is determined by $g\in K\exp(\mu(g))K$.
Its norm gives the distance from $o=eK$ to $go$ in $G/K$ for the corresponding invariant Riemannian metric. 
To retain the Cartan data associated to $\theta$, we use the \emph{partial Cartan projection}
$$
 \mu_\theta(g):=p_\theta(\mu(g))\in\mathfrak a_\theta^+.
$$
Since $\Gamma$ is $\theta$-transverse, $\mu_\theta|_\Gamma$ is proper: for every $T>0$, there are only finitely many $\gamma\in\Gamma$ with $\|\mu_\theta(\gamma)\|\leq T$.

Let $P_\theta=M_\theta A_\theta N_\theta^-$ be the corresponding standard parabolic subgroup with its Langlands decomposition, and let $\i:\mathfrak a\to\mathfrak a$ be
the opposition involution. 
Write $\mathcal F_\theta=G/P_\theta$ and $\mathcal F_{\i\theta}=G/P_{\i\theta}$. 
For each $\gamma\in\Gamma$, choose a Cartan decomposition $\gamma=k_\gamma\exp(\mu(\gamma))\ell_\gamma$, with $k_\gamma,\ell_\gamma\in K$, and define its two Cartan flags by
$$
 x_\gamma:=k_\gamma P_\theta\in\mathcal F_\theta,\qquad
 y_\gamma:=\ell_\gamma^{-1}w_0P_{\i\theta}
 \in\mathcal F_{\i\theta},
$$
where $w_0\in N_K(A)$ represents the longest Weyl element.
In rank one, when $\gamma o\ne o$, they are the endpoints at infinity of the geodesic rays from $o$ through $\gamma o$ and $\gamma^{-1}o$.
We call $x_\gamma$ and $y_\gamma$ the \emph{forward and backward Cartan flags} of $\gamma$, respectively.
These flags are independent of the chosen decomposition for all but finitely many $\gamma\in\Gamma$; for the remaining elements we fix a choice once and for all.

\subsection*{The direction of maximal growth and local mixing}
The $\theta$-growth indicator $\psi_\Gamma^\theta$ records the exponential growth of $\Gamma$ in different Cartan directions (see Section \ref{ps}). 
For Euclidean ball counting, we maximize this growth rate on the unit sphere. 
Set
$$
 \delta:=\max_{\substack{v\in\mathfrak a_\theta^+,\ \|v\|=1}}
 \psi_\Gamma^\theta(v).
$$
There is a unique unit vector $u=u_\Gamma$ attaining this maximum.
We call $u$ the \emph{direction of maximal growth} for the chosen Euclidean norm. 
The linear form
$$
 \psi(v):=\delta\langle u,v\rangle
$$
is tangent to $\psi_\Gamma^\theta$ at $u$ and is
$(\Gamma,\theta)$-proper by Lemma \ref{lem.dir}.
We use this form to choose Patterson-Sullivan measures that define our Bowen-Margulis-Sullivan measure.

To formulate the local mixing hypothesis, we first describe the flow space determined by the limit sets of $\Gamma$.
Let $\Lambda_\theta$ and $\Lambda_{\i\theta}$ be the limit sets of $\Gamma$ in $\mathcal F_\theta$ and $\mathcal F_{\i\theta}$, respectively.
Denote the open $G$-orbit in $\mathcal F_\theta\times\mathcal F_{\i\theta}$ by $\mathcal F_\theta^{(2)}$.
Its elements are pairs of flags in general position; in rank one, these are pairs of distinct points at infinity.
The Hopf parametrization identifies $G/M_\theta$ with $\mathcal F_\theta^{(2)}\times\mathfrak a_\theta$.
The right $A_\theta$-action translates the $\mathfrak a_\theta$ coordinate.

For $\theta\ne\Pi$, $M_\theta$ is noncompact, and the action of $\Gamma$ on $G/M_\theta$ need not be properly discontinuous.
Restricting the flag coordinates to the limit sets gives a subspace on which the action is properly discontinuous.
Set
$$
 \Lambda_\theta^{(2)}
 :=(\Lambda_\theta\times\Lambda_{\i\theta})
   \cap\mathcal F_\theta^{(2)}.
$$
Since $\Gamma$ is $\theta$-transverse, its action on $\Lambda_\theta^{(2)}\times\mathfrak a_\theta$ is properly discontinuous \cite[Theorem 9.1]{KOW}. 
Thus
$$
 \Omega:=\Gamma\backslash
 (\Lambda_\theta^{(2)}\times\mathfrak a_\theta)
$$
is a locally compact Hausdorff space with a right $A_\theta$-action.
 
We have $\ker\psi=u^\perp$. Let $\nu=\nu_\psi$ and $\nu_{\mathrm i}=\nu_{\psi\circ\i}$ be $(\Gamma,\psi)$- and $(\Gamma,\psi\circ\i)$-Patterson-Sullivan probability measures on $\Lambda_\theta$ and $\Lambda_{\i\theta}$, respectively.
Denote by $m$ the associated Bowen-Margulis-Sullivan (BMS) measure on $\Omega$, defined using Euclidean Lebesgue measure on $\mathfrak a_\theta$.
When $\dim\mathfrak a_\theta>1$, this measure has infinite total mass.
Local mixing describes the asymptotics of its compactly supported correlations after normalization.

We assume \emph{local mixing in the direction $u$}, with transverse perturbations of order $\sqrt t$, in the following sense. 
There exist a function $\Psi:(0,\infty)\to(0,\infty)$ and a bounded Borel function $J_u:\ker\psi\to(0,\infty)$ such that, for every fixed $w\in\ker\psi$ and $f_1,f_2\in C_c(\Omega)$,
\begin{equation*}
 \lim_{t\to\infty}\Psi(t)\int_\Omega f_1(x\exp(tu+\sqrt t\,w))f_2(x)\,dm(x)=J_u(w)m(f_1)m(f_2).
\end{equation*}
The function $J_u$ describes the dependence of the limit on the transverse parameter $w$.
The case $w=0$ concerns the directional flow $x\mapsto x\exp(tu)$. 
For fixed $w\in u^\perp$, the vectors $tu+\sqrt t\,w$ have transverse displacement of order $\sqrt t$ and, after division by $t$, converge to $u$.
For each pair $f_1,f_2$, we also require these normalized correlations to be bounded in absolute value uniformly over $t\geq1$ and $w\in\ker\psi$ with $\exp(tu+\sqrt t\,w)\in A_\theta^+$.
This local mixing hypothesis implies $u\in\operatorname{int}\mathcal L_\theta\subset\operatorname{int}\mathfrak a_\theta^+$ by Lemma \ref{lem.lm-interior}; here $\scrL_{\theta}$ is the limit cone of $\Gamma$.
 
\subsection*{Local limit theorems}
For $\gamma\in\Gamma$, write
$$
 \mu_\theta(\gamma)=t_\gamma u+z_\gamma,\qquad
 t_\gamma=\langle\mu_\theta(\gamma),u\rangle,\qquad z_\gamma\in u^\perp.
$$

Set
$$ q:=\dim\mathfrak a_\theta,\qquad r:=\frac{q-1}{2}.$$

Both local limit theorems concern elements for which $t_\gamma-T$ remains bounded.
The first describes the distribution of $z_\gamma/\sqrt T$, while the second describes the distribution of $z_\gamma$ in a fixed bounded set.
Both also give the joint distribution of the two Cartan flags.

\begin{theorem}\label{thm.lm}
Under the local mixing hypothesis, for all $\kappa_1\in C(\mathcal F_\theta)$, $\kappa_2\in C(\mathcal F_{\i\theta})$ and $F\in C_c(\mathbb R\times u^\perp)$, the following hold.
\begin{enumerate}
\item
\begin{multline*}
 \lim_{T\to\infty}\frac{\Psi(T)e^{-\delta T}}{T^r}
 \sum_{\gamma\in\Gamma}\kappa_1(x_\gamma)\kappa_2(y_\gamma)
 F\left(t_\gamma-T,\frac{z_\gamma}{\sqrt T}\right)\\
 =\nu(\kappa_1)\nu_{\mathrm i}(\kappa_2)
 \int_{\mathbb R\times u^\perp}
 e^{\delta s}F(s,w)J_u(w)\,ds\,dw.
\end{multline*}
\item
\begin{multline*}
 \lim_{T\to\infty}\Psi(T)e^{-\delta T}
 \sum_{\gamma\in\Gamma}\kappa_1(x_\gamma)\kappa_2(y_\gamma)
 F(t_\gamma-T,z_\gamma)\\
 =J_u(0)\nu(\kappa_1)\nu_{\mathrm i}(\kappa_2)
 \int_{\mathbb R\times u^\perp}
 e^{\delta s}F(s,w)\,ds\,dw.
\end{multline*}
\end{enumerate}
\end{theorem}
Theorem \ref{thm.llt} proves these statements more generally, for the directions $u$ and tangent forms $\psi$ considered in Section \ref{sec.set}.
 
\begin{Remark} \rm 
In rank one, the two local limit statements coincide. For geometrically finite groups, the BMS measure is finite (\cite[Theorem 3(i)]{Sul}, \cite[Theorem 4.1]{CI}) and strongly mixing \cite{Bab}, and Theorem \ref{thm.lm} follows from Roblin's orbital equidistribution theorem \cite[Theorem 4.1.1]{Ro}.
For infinite-index normal subgroups of convex cocompact groups with abelian quotient, the BMS measure is infinite. Oh-Pan proved local mixing with $\Psi(T)=T^{d/2}$, where $d$ is the rank of the abelian quotient, and obtained corresponding counting asymptotics \cite[Theorems 4.7 and 1.13]{OP}.
 Taking differences of their ball counts, with the factors
 $\kappa_1(x_\gamma)\kappa_2(y_\gamma)$, and approximating continuous test functions gives
 Theorem \ref{thm.lm}.
\end{Remark}
\begin{Remark}\rm
For lattices in $G$, analogues of Theorem \ref{thm.lm} for the full Cartan projection ($\theta=\Pi$) follow from the mixing and counting method of Eskin-McMullen \cite{EM} and its refinement for Cartan bisectors by Gorodnik-Oh \cite{GO}, together with the Cartan integration formula.
In this case, $\Psi\equiv1$ and $J_u$ is a positive constant, while $\nu$ and $\nu_{\mathrm i}$ are the $K$-invariant probability measures on the full flag variety.
\end{Remark}

\subsection*{Counting in balls and tubes}
We now integrate the local limits to count group elements in Euclidean balls and in their intersections with fixed-width tubes around $\mathbb Ru$.

Set $c_u:=\frac{1}{\delta}\int_{u^\perp} J_u(w)e^{-\delta\|w\|^2/2}\,dw>0$, and let $\omega_{q-1}$ denote the volume of the Euclidean unit ball in $u^\perp$, with $\omega_0=1$.
For $\epsilon>0$, let $\mathbb T_\epsilon:=\{v\in\mathfrak a_\theta: \operatorname{dist}(v,\mathbb Ru)<\epsilon\}$ be the tube of radius $\epsilon$ around the line $\mathbb Ru$.
 
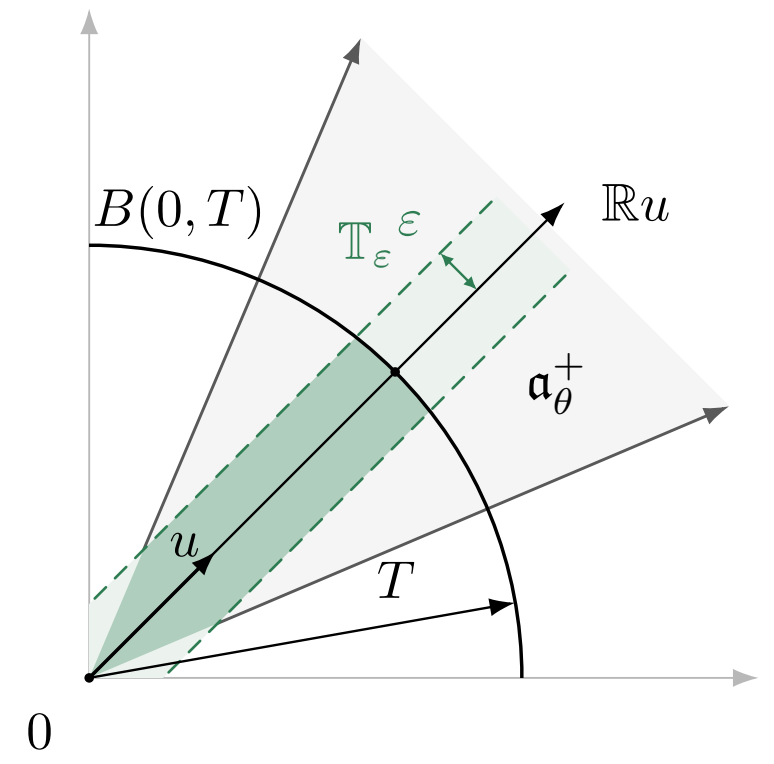
\begin{figure}[ht]
\centering
\definecolor{tubegreen}{RGB}{46,125,82}
\begin{tikzpicture}[
  x=1.08cm,y=1.08cm,
  font=\small,
  line cap=round,line join=round,
  >={Latex[length=1.7mm,width=1.2mm]}
]
\def\ballradius{2.65}
\def\tuberadius{0.32}
 
\draw[->,black!28,line width=0.35pt] (0,0)--(4.1,0);
\draw[->,black!28,line width=0.35pt] (0,0)--(0,4.1);
 
\begin{scope}
\clip (0,0) rectangle (4.25,4.25);
 
\begin{scope}[rotate=45]
 
\fill[black!4]
  (0,0)--(3.95,1.596)--(3.95,-1.596)--cycle;
\draw[->,black!65,line width=0.55pt]
  (0,0)--(3.95,1.596);
\draw[->,black!65,line width=0.55pt]
  (0,0)--(3.95,-1.596);
\node at (3.3,-0.75) {$\mathfrak a_\theta^+$};
 
\fill[tubegreen!9]
  (-3.85,-\tuberadius) rectangle (3.85,\tuberadius);
\begin{scope}
  \clip (0,0) circle[radius=\ballradius];
  \clip (0,0)--(4,1.616)--(4,-1.616)--cycle;
  \fill[tubegreen!38]
    (-4,-\tuberadius) rectangle (4,\tuberadius);
\end{scope}
 
\draw[tubegreen,line width=0.5pt,
      dash pattern=on 3pt off 2.4pt]
  (-3.85,\tuberadius)--(3.85,\tuberadius);
\draw[tubegreen,line width=0.5pt,
      dash pattern=on 3pt off 2.4pt]
  (-3.85,-\tuberadius)--(3.85,-\tuberadius);
\draw[line width=0.65pt]
  (0,0) circle[radius=\ballradius];
 
\draw[->,line width=0.45pt]
  (0,0)--(4.12,0) node[right=3pt] {$\mathbb Ru$};
\draw[->,line width=0.65pt]
  (0,0)--(1.1,0) node[pos=0.75,above=1pt] {$u$};
 
\draw[->,line width=0.45pt]
  (0,0)--(-35:\ballradius)
  node[pos=0.72,above=1pt] {$T$};
\fill (\ballradius,0) circle[radius=0.9pt];
\node[right=2pt] at (\ballradius,-0.55) {};
\draw[<->, >={Latex[length=0.9mm,width=0.8mm]},
      tubegreen,line width=0.45pt]
  (3.36,0)--(3.36,\tuberadius)
  node[midway,above left=3pt] {{$\epsilon$}}; 
 
\node[above=4pt,text=tubegreen]
  at (2.8,0.4) {{$\mathbb T_\epsilon$}}; 
 
\end{scope}
 
\node at (0.55,2.85) {$B(0,T)$};
\end{scope}
 
\fill (0,0) circle[radius=0.9pt];
\node[below left=3pt] at (0,0) {$0$};
 \end{tikzpicture}
 \caption{A fixed-width tube around the direction of maximal growth,
intersected with a Euclidean ball.} \label{fig.tube}
\end{figure}

\begin{theorem}\label{thm.cnt1}
Under the local mixing hypothesis,
the following hold.
\begin{enumerate}
\item
$$
 \#\{\gamma\in\Gamma:\|\mu_\theta(\gamma)\|<T\}
 \sim c_u\frac{T^re^{\delta T}}{\Psi(T)}.
$$
\item For every $\epsilon>0$,
$$
 \#\{\gamma\in\Gamma:\|\mu_\theta(\gamma)\|<T,\,
       \mu_\theta(\gamma)\in\mathbb T_\epsilon\}
 \sim\frac{J_u(0)\omega_{q-1}\epsilon^{q-1}}{\delta}
       \frac{e^{\delta T}}{\Psi(T)}.
$$
 \end{enumerate}
 
\end{theorem}

The tube count has no factor $T^r$,  since $z_\gamma$ remains in a fixed bounded set. In the ball count, the main contribution  comes from $z_\gamma$ of size $\sqrt T$.
We deduce (2) from Theorem \ref{thm.tubes} by comparing $\mathbb T_\epsilon\cap B(0,T)$ with sets of the form $[0,T']u+D$, where $D\subset u^\perp$ is a closed ball.
For $\kappa_1\in C(\mathcal F_\theta)$ and
$\kappa_2\in C(\mathcal F_{\mathrm i\theta})$, the corresponding normalized sums with the factor $\kappa_1(x_\gamma)\kappa_2(y_\gamma)$ converge to the
constants in Theorem \ref{thm.cnt1}(1) and (2), multiplied by $\nu(\kappa_1)\nu_{\mathrm
i}(\kappa_2)$; see Theorems \ref{thm.eqd} and
\ref{thm.tubes}.

\subsection*{Polynomial error estimates}
The asymptotic statements above require only local mixing.
To obtain polynomial errors, we also assume linear $\theta$-regularity near $u$.
This condition gives a linear lower bound for the relevant root gaps in terms of the size of the partial Cartan projection; see Definition \ref{def.lin}.
It makes the geometric comparison errors decay exponentially in $\|\mu_\theta(\gamma)\|$ near $u$.
It is used only for the error estimates.

\begin{theorem}\label{thm.cnt2}
Suppose that local mixing has a polynomial error in the sense of Definition \ref{def.em} and that $\Gamma$ is linearly $\theta$-regular near $u$ in the sense of Definition \ref{def.lin}.
Then there exist $\eta>0$ and an integer $k\geq1$ such that the errors in \Cref{thm.lm}  are
$$
 O_D\!\left(
 T^{-\eta}\|\kappa_1\|_{C^k}\|\kappa_2\|_{C^k}\|F\|_{C^1}
 \right),
$$
uniformly for smooth $\kappa_1,\kappa_2$ and $F\in C_c^1(\mathbb R\times u^\perp)$ supported in a fixed compact set $D$.
Moreover, the ball count in Theorem \ref{thm.cnt1}(1) has relative error $O(T^{-\eta})$. For each fixed $\e >0$, the tube count in Theorem \ref{thm.cnt1}(2) has relative
error $O_\e(T^{-\eta})$.
\end{theorem}

\subsection*{Applications to Anosov groups}
For $\theta$-Anosov groups, the local mixing hypothesis holds with $\Psi(t)=t^r$. 
Thirion \cite{ThiMix} proved local mixing for Ping-Pong groups along the maximal growth direction.
For Anosov groups, local mixing was proved
by Sambarino \cite{Samb1,Samb4} and Chow-Sarkar \cite{CS}; see also Edwards-Lee-Oh \cite[Theorem 7.4]{ELO} for
the Borel case.

For torsion-free $\theta$-Anosov groups, the precise statement is \cite[Proposition
A.1]{ELO2}. For general $\theta$-Anosov $\Gamma$, it follows by applying this result to a torsion-free
subgroup of finite index (Selberg's lemma) and lifting the test functions. Since
$u\in\operatorname{int}\mathcal L_\theta$ by Lemma \ref{lem.dir}, Theorem \ref{thm.cnt1} gives the following.
 
\begin{theorem}
If $\Gamma$ is $\theta$-Anosov, then, for every fixed $\epsilon>0$, as $T\to\infty$,
$$
\#\{\gamma\in\Gamma:\|\mu_\theta(\gamma)\|<T\}\sim c_u e^{\delta T};$$
$$ \#\{\gamma\in\Gamma:\|\mu_\theta(\gamma)\|<T, \mu_\theta(\gamma)\in\mathbb T_{\epsilon}\} 
 \sim\frac{J_u(0)\omega_{q-1}\epsilon^{q-1}}{\delta}
 T^{-r}e^{\delta T}. $$
\end{theorem} 
For Borel Anosov groups, ball-counting results were obtained by Sambarino \cite{Samb1} and Edwards-Lee-Oh \cite{ELO}, and tube counting by Chow-Oh \cite[Theorem 6.1]{CO}. Theorem \ref{thm.cnt1} gives both asymptotics for partial Cartan projections under a common local mixing hypothesis, with a general normalizing function $\Psi$.

\subsection*{Applications to relatively Anosov surface groups}
Kim-Oh-Pan \cite[Theorems 1.1(2) and 17.1(2)]{KOP} proved BMS local mixing for the following classes of representations, with $\Psi(t)=t^r$ and polynomial error of order $1/2$.
 
Let $\Gamma_0<\operatorname{PSL}_2(\mathbb R)$ be a torsion-free, non-elementary, finitely generated discrete subgroup with parabolic elements, and let $\rho:\Gamma_0\to G$ be a Zariski dense representation. Assume that one of the following holds:
\begin{enumerate}
\item $G$ is $\mathbb R$-split and $\rho$ is type-preserving and positive.
\item $G$ is $\mathbb R$-split and simple, $\operatorname{rank}_{\mathbb R}G\leq3$, and $\rho$ is relatively Borel Anosov.
\item $G=\prod_{j=1}^{\ell}\operatorname{SO}(n_j,1)^\circ$, with $n_j\geq2$, and $\rho$ is type-preserving with geometrically finite projections to all factors.
\end{enumerate}
Type preservation and positivity are as in \cite[Definitions 3.2 and 4.4]{KOP}. These representations are relatively Borel Anosov by \cite[Theorem 4.6 and Proposition 4.9]{KOP}. Set $\Gamma=\rho(\Gamma_0)$ and $\theta=\Pi$, and use the Euclidean norm on $\mathfrak a$ induced by the Killing form.
For these groups, the direction of maximal growth $u$ belongs to $\operatorname{int}\mathcal L_\Pi$ \cite[Section 18.4]{KOP}, so the cited local mixing theorem applies at $u$.
\begin{corollary}\label{cor.KOP}
There exist $c>0$ and $\eta>0$ such that
\begin{equation*}
 \#\{\gamma\in\Gamma:\|\mu(\gamma)\|<T\}
 =c e^{\delta T}\bigl(1+O(T^{-\eta})\bigr).
\end{equation*}
Moreover, for every $\epsilon>0$, there exists $c_\epsilon>0$ such that 
\begin{multline*}
 \#\{\gamma\in\Gamma:\|\mu(\gamma)\|<T,\ \operatorname{dist}(\mu(\gamma),\mathbb Ru)<\epsilon \}=c_\epsilon\frac{e^{\delta T}}{T^r} \bigl(1+O_\epsilon(T^{-\eta})\bigr).
\end{multline*}
\end{corollary}

\subsection*{Counting for the sup norm}
The  power of $T$
in the ball-counting asymptotic depends on the norm. For $G=\operatorname{SL}_3(\mathbb R)$, we identify $\mathfrak a$ with the space of
 trace-zero diagonal matrices and consider the sup-norm on $\fa$:
$$
 \|\operatorname{diag}(v_1,v_2,v_3)\|_\infty:=\max_i|v_i|.
$$
 
\begin{corollary}\label{cor.sup}
Suppose that $G=\operatorname{SL}_3(\mathbb R)$ and $\Gamma$ is one of the groups in Corollary \ref{cor.KOP}. There exist $\delta_\infty, c_{\infty}>0$ such that
\begin{equation}\label{eq.sup}
 \#\{\gamma\in\Gamma:\|\mu(\gamma)\|_\infty<T\}
 \sim c_{\infty}\frac{e^{\delta_\infty T}}{\sqrt T}.
\end{equation}
\end{corollary}
 
\subsection*{Related work}
In higher rank, Quint \cite{Quint3} and Thirion \cite{Thi} proved counting results for Schottky groups and Ping-Pong groups, respectively. 
Sambarino \cite{Samb1} obtained orbital counting for hyperconvex representations, building on ideas of Roblin. 
Canary, Zhang, Zhu and Zimmer \cite[Theorems 1.1 and 1.2]{CZZZ} proved counting results for elements of relatively $\theta$-Anosov groups satisfying $\varphi(\mu_\theta(\gamma))\leq T$, where $\varphi$ is a linear form with finite critical exponent.
Their proof uses a finite BMS measure on the flow space associated with $\varphi$ and requires no local mixing assumption for the full $\mathfrak a_\theta$-action. 
Their argument relies on the linearity of $\varphi$ and does not directly give counting asymptotics in Euclidean balls.

\subsection*{Proof strategy and organization}
We compare BMS correlations with smoothed counting sums, using bounds for the Busemann map and the Cartan flags. 
We then integrate against $e^{\psi(v)}\,dv$ and pass from the smoothed sums to the counting sums to prove Theorem \ref{thm.llt}. 
Counting in balls and tubes follows by integrating the local limits and controlling the tails. 
Under linear $\theta$-regularity near $u$, quantitative versions of these arguments give the stated polynomial errors.

Section \ref{ps} recalls transverse subgroups and Patterson-Sullivan measures. Section \ref{sec.cartan-comparison}  compares Cartan projections, Cartan flags and Busemann maps.
Section \ref{sec.set} states the local mixing hypotheses, proves the two local limit theorems, and deduces counting in tubes.
Section \ref{sec.eucl} proves counting in Euclidean balls and Theorem \ref{thm.cnt2}.
Section \ref{sec.sup} proves sup-norm counting and treats the $\operatorname{SL}_3(\mathbb R)$ example.

\subsection*{AI statement}
We used AI tools to assist with proofreading, reference searches, checking proofs, and drawing Figures 1-2.

\section{Transverse subgroups and Patterson-Sullivan measures}\label{ps}
 
 Let $G$ be a connected semisimple real algebraic group. Fix a maximal compact subgroup $K<G$, a compatible maximal split torus $A<G$, and a closed positive Weyl chamber $\fa^+\subset\fa:=\log A$. Let $\Pi$ be the corresponding set of simple roots. We fix a Weyl-group invariant inner product $\langle\cdot,\cdot\rangle$ on $\fa$, with norm $\|\cdot\|$, and write $o=eK\in G/K$. Let $d$ be the corresponding $G$-invariant Riemannian distance on $G/K$.
 
Fix a nonempty subset $\theta\subset\Pi$, and set $$\fa_\theta:=\bigcap_{\alpha\in\Pi\setminus\theta}\ker\alpha,\quad \fa_\theta^+:=\fa_\theta\cap\fa^+,\quad A_\theta=\exp\fa_\theta,\quad A_\theta^+=\exp\fa_\theta^+.$$
Let $p_\theta:\mathfrak a\to\mathfrak a_\theta$
denote the orthogonal projection.
For $\ell\in\mathfrak a_\theta^*$, we write $\|\ell\|_*:=\sup_{\|v\|=1}|\ell(v)|$ for the dual norm.
 
\begin{Def} For each $g\in G$, let $\mu(g)\in \mathfrak a^+$ be the unique element such that $g\in
K\exp(\mu(g))K$ \cite[Theorem 7.39]{Knapp}.
The \textit{partial Cartan projection} $\mu_\te:G\to \fa_\te^+$ is defined as
$$
 \mu_\te(g)=p_\te (\mu(g)).
$$
 
\end{Def}

Let $w_0\in N_K(A)$ represent the longest Weyl element. The opposition involution is $\i:=-\operatorname{Ad}(w_0)|_{\fa}$. It preserves $\fa^+$ and permutes $\Pi$; we denote the image of $\theta$ by $\i\theta$.
We have $\mu(g^{-1})=\i\mu(g)$ and $\mu_{\i\theta}(g^{-1})=\i\mu_\theta(g)$.
 
Let $L_\theta=Z_G(A_\theta)$ and let $M_\theta$ be the subgroup generated by $K\cap L_\theta$ and $[L_\theta,L_\theta]$. Then $L_\theta=A_\theta M_\theta$. Let $N_\theta:=N_\theta^-$ and $N_\theta^+$ be the connected subgroups whose Lie algebras are the sums of, respectively, the positive and negative root spaces with nonzero restriction to $\fa_\theta$. We set
$$ P_\theta:=P_\te^-=L_\theta N_\theta^-,\quad
 P_\theta^+:=L_\theta N_\theta^+, \quad \F_\theta:=G/P_\theta.
$$
 
For $g\in G$, set
$$
 g^+:=gP_\theta,\qquad g^-:=gw_0P_{\i\theta}.
$$
In particular, $e^+=P_\theta$ and $e^-=w_0P_{\i\theta}$.
The open $G$-orbit
$$
 \F_\theta^{(2)}:=\{(g^+,g^-):g\in G\}
 \subset\F_\theta\times\F_{\i\theta}
$$
consists of pairs of flags in \emph{general position}.
 
\subsection*{Busemann maps and Gromov products.}
Let $\sigma:G\times\F_\Pi\to\fa$ be the Iwasawa cocycle: for $\xi=kP_\Pi$ with $k\in K$, the vector $\sigma(g,\xi)$ is determined by $gk\in K\exp(\sigma(g,\xi))N_\Pi$. For $\xi\in\F_\theta$, define
$$
 \sigma_\theta(g,\xi):=p_\theta\sigma(g,\xi_0),\qquad \beta_\xi^\theta(g,h):=\sigma_\theta(g^{-1},\xi)-\sigma_\theta(h^{-1},\xi),
$$
where $\xi_0\in\F_\Pi$ projects to $\xi$. These definitions are independent of $\xi_0$ \cite[Lemma 6.1]{Quint2}. The Busemann map is continuous and depends on $g,h$ only through $gK,hK$. We therefore also write $\beta_\xi^\theta(go,ho)$.
 The
 Busemann map satisfies the following identities \cite[Lemma 6.2]{Quint2}:
\begin{equation}\label{eq.basic0}
\beta_\xi^\te(g,h)+\beta_\xi^\te(h,q)=\beta_\xi^\te(g,q),\quad \beta_{g\xi}^\te(gh,gq)=\beta_\xi^\te(h,q).
\end{equation}
 
For $(x,y)\in\mathcal F_\theta^{(2)}$, its Gromov product is given by
\begin{equation}\label{pc.H}
 \mathcal G(x,y) =\frac{1}{2} \left(\beta_x^\theta(o,ho) +\mathrm i\beta_y^{\mathrm i\theta}(o,ho)\right)
\end{equation}
for any $h\in G$ with $(h^+,h^-)=(x,y)$. The sum is independent of $h$ and satisfies
\begin{equation}\label{eq.gp}
 2\mathcal G(gx,gy)=2\mathcal G(x,y)+\beta_{gx}^\theta(o,go) +\mathrm i\beta_{gy}^{\mathrm i\theta}(o,go).
\end{equation}
 
\subsection*{Transverse subgroups.}
Let $m_\theta$ denote the $K$-invariant probability measure on $\F_\theta$. For a discrete subgroup $\Ga<G$, its limit set $\Lambda_\theta$ consists of the flags $\xi\in\F_\theta$ for which $(\gamma_n)_*m_\theta\to\delta_\xi$ for some sequence of distinct elements $\gamma_n\in\Ga$. Equivalently,
it consists of all $\xi\in\mathcal F_\theta$ for which there is a sequence $\gamma_n\to\infty$ in $\Gamma$, with Cartan decompositions $\gamma_n=k_n\exp(\mu(\gamma_n))\ell_n$, $k_n,\ell_n\in K$, such that $\min_{\alpha\in\theta}\alpha(\mu(\gamma_n))\to+\infty$ and $k_nP_\theta\to\xi$ in $\mathcal F_\theta$ \cite[Lemma 3.5 and its proof]{Ben}.

\begin{Def} A discrete subgroup $\Ga<G$ is called
\begin{enumerate}
    \item \textit{$\te$-regular} if $\min_{\alpha\in\te}\alpha(\mu(\gamma))\to+\infty$  as $\gamma\to\infty$ in $\Ga$, and
    \item \textit{$\te$-antipodal} if any two distinct points of $\La_{\te\cup\i\te}$ are in general position.
\end{enumerate}
We say $\Ga$ is \emph{$\te$-transverse} if it is both $\te$-regular and $\te$-antipodal. \end{Def}

Throughout the paper, we assume that 
$$\text{$\Ga<G$ is a Zariski dense $\theta$-transverse subgroup.}$$
 
Since $\mu(\gamma^{-1})=\i\mu(\gamma)$, the group $\Gamma$ is also $\i\theta$-regular and
hence $(\theta\cup\i\theta)$-regular.
The $\theta$-limit cone of $\Ga$ is the asymptotic cone of $\mu_\theta(\Ga):$
$$
 \mathcal L_\theta:=\left\{\lim_{n\to\infty}t_n\mu_\theta(\gamma_n)\in \fa_\te^+:\gamma_n\in\Ga,\ t_n>0,\ t_n\to0\right\}.
$$
 
Since $\Gamma$ is $\theta$-regular, $\mu_\theta|_\Gamma$ is proper. The limit cone $\mathcal L_\theta$ is convex with nonempty interior, and $\mu_\theta(\Gamma)$ is contained in a bounded neighborhood of $\mathcal L_\theta$ \cite[Lemma 3.4]{KOW}.

\subsection*{Growth indicator.}
For $v\ne0$, define
\begin{equation*}
 \psi_\Ga^\theta(v):=\|v\| \inf_{\substack{\mathcal C\subset\mathfrak a_\theta\text{ open cone},\, v\in\mathcal C}}\tau_\mathcal C,
\end{equation*}
where $\tau_\mathcal C$ is the abscissa of convergence of $\sum_{\gamma\in\Ga,\,\mu_\theta(\gamma)\in\mathcal C}e^{-s\|\mu_\theta(\gamma)\|}$. Set $\psi_\Ga^\theta(0)=0$.
The function $\psi_\Gamma^\theta$ is finite on $\mathcal L_\theta$, homogeneous of degree one, upper semicontinuous and concave. It is nonnegative on $\mathcal L_\theta$, positive on $\operatorname{int}\mathcal L_\theta$, and equal to $-\infty$ outside $\mathcal L_\theta$; see \cite[Theorem 1.1.1]{Quint1} and \cite[Theorem
3.3]{KOW}.
 A linear form $\psi\in \fa_\te^*$ is called $(\Ga,\theta)$-proper if $\psi(\mu_\theta(\gamma))\to+\infty$ as $\gamma\to\infty$ in $\Ga$. If $\psi$ is positive on $\mathcal L_\theta\setminus\{0\}$, then it is $(\Ga, \theta)$-proper \cite[Lemma 4.3(1)]{KOW}.
A linear form $\psi\in\fa_\theta^*$ is called a tangent form at $u\in\mathcal L_\theta\setminus\{0\}$ if $$\psi\geq\psi_\Ga^\theta \qquad \text{and}\qquad \psi(u)=\psi_\Ga^\theta(u).$$
 
\subsection*{Patterson-Sullivan measures.} 
 For $\psi\in\fa_\theta^*$, a $(\Ga,\psi)$-Patterson-Sullivan measure is a probability measure $\nu_\psi$ supported on $\Lambda_\theta$ such that for all $\gamma\in\Ga$,
\begin{equation}\label{eq.conf}
\frac{d\gamma_*\nu_\psi}{d\nu_\psi}(\xi)=e^{\psi(\beta_\xi^\theta(e,\gamma))}
\end{equation}
where $\ga_* \nu_\psi(Q)=\nu_{\psi}(\gamma^{-1} Q)$ for any Borel subset $Q\subset \F_\theta$.
By \eqref{eq.conf}, $\operatorname{supp}\nu_\psi$ is a nonempty closed
$\Gamma$-invariant subset of $\Lambda_\theta$.
Since $\Lambda_\theta$ is $\Gamma$-minimal \cite[Lemma 5.2(3)]{KOW},
we have $\operatorname{supp}\nu_\psi=\Lambda_\theta$.
If $\psi$ is a $(\Gamma,\theta)$-proper tangent form at $u\in\mathcal L_\theta\setminus\{0\}$ with $\psi_\Gamma^\theta(u)>0$, then a $(\Ga,\psi)$-PS measure exists by \cite[Proposition 5.12 and the discussion following it]{KOW} (see also \cite{CZZ}).

\subsection*{Relatively $\theta$-Anosov subgroups}
A $\theta$-transverse group $\Gamma$ is called \textit{relatively $\theta$-Anosov} if every $\xi\in\Lambda_\theta$ satisfies one of the following:
\begin{enumerate}
\item $\xi$ is conical: there are distinct $a,b\in\Lambda_\theta$ and a sequence $\gamma_n\to\infty$ in $\Gamma$ such that $\gamma_n\xi\to a$ and $\gamma_n\zeta\to b$ uniformly for $\zeta$ in compact subsets of $\Lambda_\theta\setminus\{\xi\}$.
\item $\xi$ is bounded parabolic: its stabilizer is infinite and acts properly on $\Lambda_\theta\setminus\{\xi\}$ with compact quotient. We also require this stabilizer
to be finitely generated.
\end{enumerate}
See \cite[Definition 7.7, Theorem 7.10]{KL}.
The \textit{$\theta$-Anosov} case is the special case in which every point of $\Lambda_\theta$ satisfies (1) \cite[Section 3.3, Theorem 1.1]{KLP}.

\begin{Def}\label{def.lin}
For $u\in\mathfrak a_\theta^+\setminus\{0\}$, we say that $\Gamma$ is \textit{linearly $\theta$-regular near $u$} if there are an open cone $\mathcal U\subset\mathfrak a_\theta$ containing $u$ and constants $c,C>0$ such that for every $\gamma\in\Gamma$ with $\mu_\theta(\gamma)\in\mathcal U$,
\begin{equation}\label{eq.rct}
 \min_{\alpha\in\theta}\alpha(\mu(\gamma))\geq c\|\mu_\theta(\gamma)\|-C.
\end{equation}
\end{Def}
This condition is used only  to obtain
polynomial errors.
If $\Gamma$ is contained in a relatively $\theta$-Anosov subgroup, then it holds for every $u\in\mathfrak a_\theta^+\setminus\{0\}$ \cite[Theorem 1.7, Propositions 5.1 and 13.4]{ZZ1}. 
When $\theta=\Pi$, every discrete subgroup $\Gamma<G$ is linearly $\Pi$-regular near every $u\in\operatorname{int}\mathfrak a^+$.
\section{Cartan comparison and directional flows}\label{sec.cartan-comparison}
We  describe the flow space and  show that every topologically transitive
direction belongs to $\operatorname{int}\mathcal L_\theta$. We then  compare Cartan projections and Busemann maps, with
bounds used in the local limit  theorems.
 
\subsection*{Hopf coordinates and the flow space.}
We consider the following Hopf parametrization:
\begin{equation}\label{eq.hopf0}
 G/M_\theta\to\F_\theta^{(2)}\times\fa_\theta,\qquad gM_\theta\mapsto(g^+,g^-,\beta_{g^+}^\theta(e,g)).
\end{equation}
 For $a\in A$ and $m\in M_\theta$,
 we have
\begin{equation*}
\begin{aligned}
 \beta_{g^+}^\theta(e,ga m)=\beta_{g^+}^\theta(e,g)+p_\theta\log a,\qquad
 \beta_{g^-}^{\i\theta}(e,ga m)=\beta_{g^-}^{\i\theta}(e,g)-\i p_\theta\log a.
 \end{aligned}
\end{equation*}
In particular, each expression with $a=e$ is independent of the representative of $gM_\theta$.
Set
$$
 \Lambda_\theta^{(2)}:=(\Lambda_\theta\times\Lambda_{\i\theta})\cap\F_\theta^{(2)},\qquad \widetilde\Omega:=\Lambda_\theta^{(2)}\times\fa_\theta.
$$
 
In Hopf coordinates, the left $\Gamma$-action and the right $A_\te$-action are given, for $g\in\Gamma$ and $v\in\mathfrak a_\te$, by
\begin{equation}\label{eq.hopfa}
 g(\xi,\eta,b) =(g \xi,g \eta,b+\beta_{g \xi}^{\te}(e,g)),\qquad (\xi,\eta,b)\exp v=(\xi,\eta,b+v).
\end{equation}
 
 The action of $\Gamma$ on $\widetilde\Omega$ is  properly discontinuous by \cite[Theorem 9.1]{KOW}, and its stabilizers are finite. We set $$\Omega:=\Gamma\backslash\widetilde\Omega$$ and denote the quotient map $\widetilde \Omega\to \Omega$ by $\pi$.
In particular, $\Omega$ is locally compact and Hausdorff.
 
\subsection*{Topological transitivity.}
The
proof of \cite[Lemma 3.2]{ELO2} gives the following bound; it does not use
torsion-freeness.
\begin{lemma}\label{lem.cpt}
For every compact set $Q\subset\widetilde\Omega$, there exists $C_Q>0$ such that, for all $\gamma\in\Gamma$, $v\in\mathfrak a_\theta$ and $z,z'\in Q$ satisfying $\gamma z=z'\exp(v)$, we have
$$
 \min\left\{ \|\mu_\theta(\gamma)-v\|, \|\mu_\theta(\gamma^{-1})+v\| \right\}\leq C_Q.
$$
\end{lemma}
For $u\in \fa_\theta^+$, we say that the flow $x\mapsto x\exp(tu)$ is \emph{topologically transitive} on $\Omega$ if, for every pair of nonempty open subsets $U,V\subset\Omega$, there exists $t\geq0$ such that $U\exp(tu)\cap V\neq\emptyset$.
\begin{lemma}\label{lem.lm-interior}
Let $u\in\mathfrak a_\theta^+$ be a nonzero vector.
If the flow $x\mapsto x\exp(tu)$ is topologically transitive on $\Omega$, then $u\in\operatorname{int}\mathcal L_\theta$.
In particular, $u\in\operatorname{int}\mathfrak a_\theta^+$.
\end{lemma}
 
\begin{proof}
Suppose that $u\notin\operatorname{int}\mathcal L_\theta$.
Since $\mathcal L_\theta$ is a closed convex cone with nonempty interior, the separating hyperplane theorem gives a nonzero linear form $\ell\in\mathfrak a_\theta^*$ such that
$$
 \ell\geq0\quad\text{on }\mathcal L_\theta, \quad\text{ and }\quad \ell(u)\leq0.
$$
By \cite[Lemma 3.4]{KOW}, $\mu_\theta(\Gamma)$ lies within a bounded distance of $\mathcal L_\theta$.
Hence there is $B>0$ such that for all $\ga\in\Ga$,
$$
 \ell(\mu_\theta(\gamma))\geq-B.
$$
Choose a nonempty relatively compact open subset $U\subset\Omega$ and a compact set $Q\subset\widetilde\Omega$ such that $U\subset\pi(Q)$.
Let $C_Q$ be the constant in Lemma \ref{lem.cpt}, and choose $v\in\mathfrak a_\theta$ such that
$$
 \ell(v)<-B-C_Q\|\ell\|_*-1.
$$
By topological transitivity, for every $T>0$, there exists $t\geq T$ such that
$$
 U\exp(tu+v)\cap U\ne\emptyset.
$$
Hence there exist $t_n\to+\infty$, $\gamma_n\in\Gamma$, and $z_n,z_n'\in Q$ such that $\gamma_nz_n=z_n'\exp(t_nu+v)$.
Lemma \ref{lem.cpt} gives
$$
 \min\left\{ \|\mu_\theta(\gamma_n)-(t_nu+v)\|, \|\mu_\theta(\gamma_n^{-1})+(t_nu+v)\| \right\}\leq C_Q.
$$
If the second term were at most $C_Q$ for infinitely many $n$, then, along a subsequence,
$$
 \frac{\mu_\theta(\gamma_n^{-1})}{t_n}\to -u.
$$
Since $\mathfrak a_\theta^+$ is closed, this would imply $-u\in\mathfrak a_\theta^+$, contradicting $0\ne u\in\mathfrak a_\theta^+$.
Thus, for all sufficiently large $n$,
$$
 \|\mu_\theta(\gamma_n)-(t_nu+v)\|\leq C_Q.
$$
Since $\ell(u)\leq0$, it follows that
$$
 \ell(\mu_\theta(\gamma_n)) \leq t_n\ell(u)+\ell(v)+C_Q\|\ell\|_*<-B,
$$
a contradiction. Thus $u\in\operatorname{int}\mathcal L_\theta$.
\end{proof}

\begin{Remark}\rm
For $\theta=\Pi$ and $u\in\operatorname{int}\mathfrak a^+$,
Dang-Glorieux \cite[Proposition 4.7]{DG} proved the implication in Lemma \ref{lem.lm-interior} for every Zariski dense discrete subgroup $\Gamma<G$.
The same argument works for $\theta=\Pi$ without the transverse assumption.
By \cite[Sections 1.2 and 4.6]{Ben}, $\mu(\Gamma)$ is contained in a bounded neighborhood of $\mathcal L_\Pi$.
Thus a linear form $\ell$ that is nonnegative on $\mathcal L_\Pi$ satisfies $\ell(\mu(\gamma))\geq-B$ for all $\gamma\in\Gamma$, for some $B>0$.
Since $M_\Pi$ is compact, $\Gamma$ acts properly on $\widetilde\Omega$.
Lifting $Q$ to a compact subset of $G$ gives
$$
 \gamma_n=g_n'\exp(t_nu+v)m_ng_n^{-1}, \quad m_n\in M_\Pi, \qquad \|\mu(\gamma_n)-(t_nu+v)\|\leq C
$$
for all sufficiently large $n$, since $t_nu+v\in\mathfrak a^+$.
 Choosing $\ell(v)<-B-C\|\ell\|_*-1$ gives the same contradiction.
\end{Remark}
 
\subsection*{Quantitative Cartan comparison.}
 
Fix $K$-invariant Riemannian distances $d_{\mathcal F_\theta}$ and $d_{\mathcal F_{\mathrm i\theta}}$ on $\cal F_\te$ and $\cal F_{\i\te}$.
For each $g\in G$, choose a Cartan decomposition $g=k_g\exp(\mu(g))\ell_g\in KA^+K$, and define the \emph{forward and backward Cartan flags} of $\gamma$ as in the introduction:
\begin{equation}\label{eq.cf1}
 x_g:=k_gP_\theta,\qquad y_g:=\ell_g^{-1}w_0P_{\mathrm i\theta}.
\end{equation}
  The flags $x_g$ and $y_g$ are independent of the choice if $\alpha(\mu(g))>0$ for all $\alpha\in\theta$.
By $\theta$-regularity, these flags are well
defined for all but finitely many $\gamma\in\Gamma$.
Moreover, as $\gamma\to\infty$ in $\Gamma$, every accumulation point of $x_\gamma$ lies in $\Lambda_\theta$, and every accumulation point of $y_\gamma$ lies in $\Lambda_{\mathrm i\theta}$. See \cite[Lemma 3.5]{Ben}.
 
For $g\in G$ we write
$$\gap(g):=\min_{\alpha\in\theta}\alpha(\mu(g))\geq0 .
$$
By $\theta$-regularity, $\gap(\gamma)\to\infty$ as $\gamma\to\infty$ in $\Gamma$.
 
\begin{lemma}\label{lem.dist1}
Let $Q\subset\mathcal F_\theta^{(2)}$ be a compact subset. There is $C_Q>0$ such that for all $g\in
G$ and  $x\in\mathcal F_\theta$ such that $(x,y_g)\in Q$,
\begin{equation}\label{eq.dist1}
 \begin{gathered} d_{\mathcal F_\theta}(gx,x_g)\leq C_Qe^{-\gap(g)} \quad\text{and}\\
 \bigl\| \mu_\theta(g)-\beta_{gx}^\theta(o,go)-2\mathcal G(x,y_g)\bigr\|\leq C_Qe^{-\gap(g)}. \end{gathered}
\end{equation}
\end{lemma}
\begin{proof}
Write $g=ka\ell$ with $k=k_g$, $\ell=\ell_g$ and $a=\exp\mu(g)$.
Suppose that $(x,y_g)\in Q$.
Since $\ell y_g=w_0P_{\mathrm i\theta}=e^-$, we have $(\ell x,e^-)\in KQ$, which is a compact subset of $\mathcal F_\theta^{(2)}$.
The set of flags in $\mathcal F_\theta$ in general position with $e^-$ is the open cell $N_\theta^+e^+$, and $n\mapsto ne^+$ is a diffeomorphism of $N_\theta^+$
onto it.
Hence $\ell x=ne^+$ for a unique $n\in N_\theta^+$, and $n$ lies in a compact subset of $N_\theta^+$ depending only on $Q$.
Set $n_a:=ana^{-1}$.
Now $\log N_\theta^+$ is the sum of the root spaces $\mathfrak g_{-\beta}$ over all positive roots $\beta$ whose restriction to $\mathfrak a_\theta$ is nonzero; each such $\beta$ is a nonnegative integral combination of simple roots with a positive coefficient on some $\alpha\in\theta$, so that $\beta(\mu(g))\geq\gap(g)$.
 Conjugation by $a$ multiplies the $\mathfrak g_{-\beta}$-component of $\log n$ by $e^{-\beta(\mu(g))}$. Since $n$ lies in a fixed compact set, we obtain
\begin{equation}\label{eq.dist2}
 d(o,n_ao)+d_{\mathcal F_\theta}(n_ae^+,e^+)\leq C_Qe^{-\gap(g)}.
\end{equation}
Since $gx=kane^+=kn_ae^+$ and $x_g=ke^+$, and $K$ acts by isometries on $\mathcal F_\theta$, the first bound in \eqref{eq.dist1} follows.
 
For the second bound, note that $n$ fixes $e^-$ and that $\beta^{\mathrm i\theta}_{e^-}(o,no)=0$, since $n^{-1}w_0\in w_0N^-_{\mathrm i\theta}\subset KN_{\mathrm i\theta}^-$.
 Since $\ell\in K$ fixes $o$, \eqref{eq.gp} and \eqref{pc.H} with $h=n$  give
$$
 2\mathcal G(x,y_g)=2\mathcal G(\ell x,\ell y_g)=2\mathcal G(ne^+,e^-)=\beta_{ne^+}^\theta(o,no)=\beta_{n_ae^+}^\theta(ao,n_aao),
$$
where the last equality is \eqref{eq.basic0} applied with $a$, using $ane^+=n_ae^+$ and $ano=n_aao$.
Similarly, $$\beta_{gx}^\theta(o,go)=\beta_{kn_ae^+}^\theta(ko,kao)=\beta_{n_ae^+}^\theta(o,ao).$$
Adding, and using the cocycle identity together with $\beta_{n_ae^+}^\theta(n_ao,n_aao)=\beta_{e^+}^\theta(o,ao)=\mu_\theta(g)$ (by \eqref{eq.basic0} and the definition of $\beta^\theta$, since $a^{-1}\in K\exp(-\mu(g))N_\Pi$), we obtain
$$
 \beta_{gx}^\theta(o,go)+2\mathcal G(x,y_g)=\beta_{n_ae^+}^\theta(o,n_aao)=\beta_{n_ae^+}^\theta(o,n_ao)+\mu_\theta(g).
$$
Since $\|\beta_\xi^\theta(p,q)\|\leq d(p,q)$ for all $\xi\in\mathcal F_\theta$ and $p,q\in G/K$, the second bound in \eqref{eq.dist1} follows from \eqref{eq.dist2}.
\end{proof}
 
The natural projections of $\Lambda_{\theta\cup\i\theta}$ onto $\Lambda_\theta$ and $\Lambda_{\i\theta}$ are homeomorphisms \cite[Lemma 9.5]{KOW}.
Write $\xi^*\in\Lambda_{\i\theta}$ for the flag corresponding to $\xi\in\Lambda_\theta$.
 
\begin{lemma}\label{lem.dist2}
Let $\psi\in\mathfrak a_\theta^*$ be a $(\Gamma,\theta)$-proper linear form.
For compact subsets $Q_-,Q_+\subset\Lambda_\theta^{(2)}$, there exist $C,c,R_0>0$ such that the following holds.
If $R\geq R_0$, $\gamma\in\Gamma$, $x\in\Lambda_\theta$, and $y\in\Lambda_{\mathrm i\theta}$ satisfy
$$
 (x,\gamma^{-1}y)\in Q_-, \quad (\gamma x,y)\in Q_+, \quad \psi\bigl(2\mathcal G(x,\gamma^{-1}y) +\beta_{\gamma x}^\theta(o,\gamma o)\bigr)\geq R,$$
then
\begin{equation}\label{eq.err1}
 \bigl\|\mu_\theta(\gamma)-2\mathcal G(x,\gamma^{-1}y)
       -\beta_{\gamma x}^\theta(o,\gamma o)\bigr\|
 \leq Ce^{-\gap(\gamma)}
\end{equation}
and
\begin{equation}\label{eq.err2}
 d_{\mathcal F_\theta}(\gamma x,x_\gamma)
 +d_{\mathcal F_{\mathrm i\theta}}(\gamma^{-1}y,y_\gamma)
 \leq Ce^{-\gap(\gamma)}.
\end{equation} Moreover, $\|\mu_\theta(\gamma)\|\geq cR$.
\end{lemma}
 
\begin{proof}
We first claim that there exist $R_0>0$ and a compact subset $Q\subset\mathcal F_\theta^{(2)}$ such that, whenever $R\geq R_0$ and $x,y,\gamma$ satisfy the hypotheses of the lemma, both $(x,y_\gamma)$ and $(x_\gamma,y)$ belong to $Q$.
 
Suppose not. Then for each $n\in\mathbb N$, there exist $\gamma_n\in\Gamma$, $x_n\in\Lambda_\theta$, and $y_n\in\Lambda_{\mathrm i\theta}$ satisfying
$$
 (x_n,\gamma_n^{-1}y_n)\in Q_-, \;\; (\gamma_nx_n,y_n)\in Q_+,\;\; \psi\bigl(2\mathcal G(x_n,\gamma_n^{-1}y_n) +\beta_{\gamma_nx_n}^\theta(o,\gamma_no)\bigr)\geq n,
$$
such that at least one of $(x_n,y_{\gamma_n})$ and $(x_{\gamma_n},y_n)$ has distance less than $1/n$ from the complement of $\mathcal F_\theta^{(2)}$.
After passing to a subsequence, the same one of these two pairs converges to a pair outside $\mathcal F_\theta^{(2)}$.
 
Since $\mathcal G(x_n,\gamma_n^{-1}y_n)$ is uniformly bounded and $\|\beta_{\gamma_nx_n}^\theta(o,\gamma_no)\|\le d(o, \gamma_n o)$, it follows that $\gamma_n\to\infty$.
Passing to a subsequence, we may assume
$$
 x_n\to x,\qquad \gamma_n^{-1}y_n\to z,\qquad \gamma_nx_n\to x',\qquad y_n\to y.
$$
Then $(x,z)\in Q_-$ and $(x',y)\in Q_+$. Write $\gamma_n=k_n\exp(\mu(\gamma_n))\ell_n\in K \exp \fa^+ K$; so
$
 \gamma_n^{-1}  =(\ell_n^{-1}w_0)\exp(\mathrm i\mu(\gamma_n)) (w_0^{-1}k_n^{-1}) \in  K \exp \fa^+ K.
$
After passing to a further subsequence, assume
$k_n\to k$ and $\ell_n\to\ell$. Then
$$
 kP_{\theta\cup\mathrm i\theta} \in\Lambda_{\theta\cup\mathrm i\theta}, \qquad \ell^{-1}w_0P_{\theta\cup\mathrm i\theta} \in\Lambda_{\theta\cup\mathrm i\theta}.
$$
Set $\xi^+=kP_\theta$ and $\xi^-=\ell^{-1}w_0P_\theta$.
The correspondence $\xi\mapsto\xi^*$ gives
$$
 (\xi^+)^*=kP_{\mathrm i\theta}, \qquad (\xi^-)^*=\ell^{-1}w_0P_{\mathrm i\theta}.
$$
 For all sufficiently large $n$, the Cartan flags are uniquely determined. Hence
\begin{align*}
 x_{\gamma_n}=k_nP_\theta&\to\xi^+,
 &
 y_{\gamma_n^{-1}}=k_nP_{\mathrm i\theta}
 &\to(\xi^+)^*,
 \\
 x_{\gamma_n^{-1}}=\ell_n^{-1}w_0P_\theta
 &\to\xi^-,
 &
 y_{\gamma_n}=\ell_n^{-1}w_0P_{\mathrm i\theta}
 &\to(\xi^-)^*.
\end{align*}
 
By the choice of the sequence, at least one of $(x,(\xi^-)^*)$ and $(\xi^+,y)$ is not in general position.
We first show that $x=\xi^-$.
Otherwise, antipodality gives $(x,(\xi^-)^*)\in\mathcal F_\theta^{(2)}$, so $(x_n,y_{\gamma_n})$ lies in a fixed compact subset of $\mathcal F_\theta^{(2)}$ for all sufficiently large $n$.
Lemma \ref{lem.dist1} and $\theta$-regularity give
$$
 d_{\mathcal F_\theta}(\gamma_nx_n,x_{\gamma_n})  \ll e^{-\gap(\gamma_n)}\to0.
$$
Hence $x'=\xi^+$. Since $(x',y)\in Q_+$, the pair $(\xi^+,y)$ is also in general position, a contradiction.
Thus $x=\xi^-$.
 
We next show that $y=(\xi^+)^*$.
Otherwise, antipodality gives $(\xi^+,y)\in\mathcal F_\theta^{(2)}$, so $(x_{\gamma_n},y_n)$ lies in a fixed compact subset of $\mathcal F_\theta^{(2)}$ for all sufficiently large $n$.
For $\mathrm i\theta$, the forward and backward Cartan flags of $\gamma_n^{-1}$ are $y_{\gamma_n}$ and $x_{\gamma_n}$, respectively, and
$$
 \min_{\alpha\in\mathrm i\theta} \alpha(\mu(\gamma_n^{-1}))=\gap(\gamma_n).
$$
Applying Lemma \ref{lem.dist1} to $\gamma_n^{-1}$ with $\theta$ replaced by $\mathrm i\theta$, we obtain
$$
 d_{\mathcal F_{\mathrm i\theta}} (\gamma_n^{-1}y_n,y_{\gamma_n}) \ll e^{-\gap(\gamma_n)}\to0.
$$
Hence $z=(\xi^-)^*$, contradicting $(x,z)\in Q_-$ because $x=\xi^-$ and $(\xi^-,(\xi^-)^*)$ is not in general position.
Therefore $y=(\xi^+)^*$.
 
It follows that
$$
 (\gamma_nx_n,y_{\gamma_n^{-1}}) \to(x',(\xi^+)^*)\in Q_+ \subset\mathcal F_\theta^{(2)}.
$$
We can now apply Lemma \ref{lem.dist1} to $g=\gamma_n^{-1}$ and the flag $\gamma_nx_n$, with $\theta$ unchanged.
Since $\gap(\gamma_n^{-1})\to\infty$, this gives
$$
 \beta_{x_n}^\theta(o,\gamma_n^{-1}o) =\mu_\theta(\gamma_n^{-1}) -2\mathcal G(\gamma_nx_n,y_{\gamma_n^{-1}})+o(1).
$$
Since both $\mathcal G(x_n,\gamma_n^{-1}y_n)$ and $\mathcal G(\gamma_nx_n,y_{\gamma_n^{-1}})$ are bounded, and $\beta_{\gamma_nx_n}^\theta(o,\gamma_no) =-\beta_{x_n}^\theta(o,\gamma_n^{-1}o)$, we obtain
$$
 2\mathcal G(x_n,\gamma_n^{-1}y_n)+\beta_{\gamma_nx_n}^\theta(o,\gamma_no)=-\mu_\theta(\gamma_n^{-1})+O(1).
$$
Since $\gamma_n^{-1}\to\infty$, the $(\Gamma,\theta)$-properness of $\psi$ implies $\psi(\mu_\theta(\gamma_n^{-1}))\to+\infty$.
 Applying $\psi$ to the preceding identity gives a contradiction.
This proves the claim.
 
Now let $R\geq R_0$, and let $x,y,\gamma$ satisfy the hypotheses of the lemma.
By the claim and Lemma \ref{lem.dist1}, applied to $\gamma$ and to $\gamma^{-1}$ with $\theta$ replaced by $\mathrm i\theta$ in the latter case,
$$
 d_{\mathcal F_\theta}(\gamma x,x_\gamma)+d_{\mathcal F_{\mathrm i\theta}}(\gamma^{-1}y,y_\gamma) \leq Ce^{-\gap(\gamma)}.
$$
This proves \eqref{eq.err2}.
The pairs $(x,\gamma^{-1}y)$ and $(x,y_\gamma)$ belong to the fixed compact set $Q_-\cup Q\subset\mathcal F_\theta^{(2)}$.
Since $\mathcal G$ is smooth, it follows that
$$
 \|\mathcal G(x,\gamma^{-1}y)-\mathcal G(x,y_\gamma)\| \leq C\,d_{\mathcal F_{\mathrm i\theta}}(\gamma^{-1}y,y_\gamma)\leq Ce^{-\gap(\gamma)}.
$$
Combining this with the second bound in \eqref{eq.dist1} proves \eqref{eq.err1}. 
Finally, \eqref{eq.err1} and the hypothesis give
$$
 \|\psi\|_*\,\|\mu_\theta(\gamma)\| \geq\psi(\mu_\theta(\gamma))\geq R-O(1).
$$
After increasing $R_0$, this gives
$\|\mu_\theta(\gamma)\|\geq cR$.
\end{proof}
\section{Local limit theorems and tube counting}\label{sec.set}
Fix $u\in\mathcal L_\theta\setminus\{0\}$ with $\psi_\Gamma^\theta(u)>0$, and a $(\Gamma,\theta)$-proper tangent form $\psi$ at $u$. Let $\nu=\nu_\psi$ and $\nu_{\mathrm i}=\nu_{\psi\circ\i}$ be PS probability measures on $\La_\theta$ and $\La_{\i\theta}$, respectively. On $\widetilde\Omega$, set
$$
 d\widetilde m(gM_\theta):=e^{\psi(\beta_{g^+}^\theta(e,g))+\psi\circ\i(\beta_{g^-}^{\i\theta}(e,g))}\,d\nu(g^+)\,d\nu_{\mathrm i}(g^-)\,db,
$$
where $db$ is Euclidean Lebesgue measure on $\mathfrak a_\theta$. This measure is left $\Gamma$-invariant and right $A_\theta$-invariant. It induces the BMS measure $m=m^{\BMS}$ on $\Omega$ such that for every $\widetilde f\in C_c(\widetilde\Omega)$,
$$
 \int_\Omega
 \left(\sum_{\gamma\in\Gamma}\widetilde f(\gamma z)\right)
 \,dm(\Gamma z)
 =\int_{\widetilde\Omega}\widetilde f(z)\,d\widetilde m(z).
$$

Set $W=\ker\psi$ and $a_{t,w}=\exp(tu+\sqrt t\,w)$ for $t>0$ and $w\in W$.
 
\begin{Def}[Local mixing]\label{lmp}\rm
We say that $m$ satisfies the local mixing property if there exist a function $\Psi:(0,\infty)\to(0,\infty)$ and a bounded positive Borel function $J_u:W\to(0,\infty)$ such that:
\begin{enumerate}
\item For every $w\in W$ and $f_1,f_2\in C_c(\Omega)$,
$$
 \lim_{t\to \infty} \Psi(t)\int_\Omega f_1(xa_{t,w})f_2(x)\,dm(x)= J_u(w)m(f_1)m(f_2).
$$
\item For every $f_1,f_2\in C_c(\Omega)$ there is $C=C(f_1,f_2)>0$ such that
$$
 \left|\Psi(t)\int_\Omega f_1(xa_{t,w})f_2(x)\,dm(x)\right|\leq C
$$
for all $t\geq1$ and $w\in W$ with $a_{t,w}\in A_\theta^+$.
\end{enumerate}
\end{Def}
 
The full support of $\nu$ and $\nu_{\mathrm i}$ and the defining
formula
for $\widetilde m$ imply that $m$ has full support on $\Omega$.
Hence local mixing along $\exp(tu)$ implies topological transitivity.  By Lemma \ref{lem.lm-interior},
$u\in\operatorname{int}\mathcal L_\theta\subset\operatorname{int}\mathfrak a_\theta^+$.
Thus $a_{t,w}\in A_\theta^+$ for every fixed $w\in W$ and all sufficiently large $t$.
 
We recall the following properties of the normalizing function $\Psi$.
The proof in \cite{ELO2} does not use that $\Gamma$ is torsion-free.
\begin{lemma}[{\cite[Lemma 6.3]{ELO2}}]\label{lem.ratio}
Suppose that $m^{\BMS}$ satisfies local mixing with normalizing function $\Psi$.
\begin{enumerate}
\item For every $C>0$, we have
$$
 \sup_{|s|\leq C} \left|\frac{\Psi(T+s)}{\Psi(T)}-1\right|\to0 \qquad \text{as $T\to\infty$.}
$$
\item Moreover, for every $\eta>0$ there is $t_\eta>0$ such that
$$
 e^{-\eta(T-t+1)} \leq\frac{\Psi(T)}{\Psi(t)} \leq e^{\eta(T-t+1)} \qquad\text{ for all $t_\eta\leq t\leq T$}.
$$
In particular, $\log\Psi(T)=o(T)$.
 \item Suppose that the local mixing property has a polynomial error of order $d>0$ (Definition \ref{def.em}). Then there exist $C_\Psi,T_0>0$ such that
$$
 \left|\log\frac{\Psi(T+s)}{\Psi(T)}\right|\leq C_\Psi(|s|+1)T^{-d} \qquad\text{for all $T\geq T_0$ and $|s|\leq T/2$.}
$$
\end{enumerate}
\end{lemma}

 Write
$$
 \err(a_{t,w},f_1,f_2):=\Psi(t)\int_\Omega f_1(xa_{t,w})f_2(x)\,dm(x)-J_u(w)m(f_1)m(f_2).
$$
 For the quantitative statements, we use smooth functions on $G/M_\theta$. If $\widetilde f\in C_c(G/M_\theta)$, define
\begin{equation}\label{eq.per}
 f(\Gamma z):=\sum_{\gamma\in\Gamma}\widetilde f(\gamma z),\qquad z\in\widetilde\Omega.
\end{equation}
 The set $\widetilde\Omega$ is closed in $G/M_\theta$, so $\operatorname{supp}\widetilde f\cap\widetilde\Omega$ is compact.
Properness of the $\Gamma$-action on $\widetilde\Omega$ makes the sum locally finite. The image in $\Omega$ of $\operatorname{supp}\widetilde f\cap\widetilde\Omega$ is compact and contains $\operatorname{supp}f$, so $f\in C_c(\Omega)$. 
 
\begin{Def}[Polynomial error]\label{def.em}\rm
The local mixing property has a  polynomial error of order $d>0$ if, for every compact set $L\subset G/M_\theta$, there are an integer $k\geq1$ and $C_L>0$ such that
\begin{equation}\label{eq.kd2}
 |\err(a_{t,w},f_1,f_2)|\leq C_Lt^{-d}\|\widetilde f_1\|_{C^k}\|\widetilde f_2\|_{C^k}
\end{equation}
whenever $t\geq1$, $a_{t,w}\in A_\theta^+$, and $\widetilde f_1,\widetilde f_2\in C_c^\infty(G/M_\theta)$ have support in $L$, with $f_1,f_2$ given by \eqref{eq.per}. The $C^k$-norms are taken in a fixed finite collection of smooth coordinate neighborhoods of $L$ in $G/M_\theta$.
\end{Def}
Lemma \ref{lem.ratio}(3) is proved as in \cite{ELO2}, applying \eqref{eq.kd2} with $w=0$ to a nonnegative $f$ as in \eqref{eq.per} with $\widetilde f$ smooth and $m(f)>0$, and to its translates $x\mapsto f(x\exp(su))$, $|s|\leq1$.

\subsection*{Local limit theorems for partial Cartan projections}
\label{sec.pc}
 Let $\|\cdot\|$ be the Euclidean norm fixed in
Section \ref{ps}, and set
$$
  \delta_u:=\psi(u)>0,\qquad
 \mathfrak j_u:=\frac{\delta_u}{\|\psi\|_*}.
$$  
Set
\begin{equation}\label{eq.qr}
q=\dim\mathfrak a_\theta=|\theta|,\qquad r=\sfrac{q-1}{2}, \qquad W=\ker\psi.    
\end{equation}
 All Lebesgue measures on $\mathfrak a_\theta$ and $W$ are induced by the Euclidean inner product. If $W=\{0\}$, its Lebesgue measure assigns mass one to the origin. 
 Since $\psi(u)> 0$, we have $\mathfrak a_\theta=\mathbb Ru\oplus W$. 
 Writing $v=tu+z$ with $z\in W$, the Euclidean volume measures satisfy
\begin{equation}\label{eq.ju}
 dv=\mathfrak j_u\,dt\,dz,\qquad \psi(v)=\psi(u) t .
\end{equation}
 
Suppose that $m^{\BMS}$ satisfies the local mixing property, and denote the functions in Definition \ref{lmp} by $(\Psi,J_u)$. 
For each $\gamma\in\Gamma$, write
$$
 \mu_\theta(\gamma)=t_\gamma u+z_\gamma,\qquad t_\gamma=\frac{\psi(\mu_\theta(\gamma))}{\psi(u)},\qquad z_\gamma\in W.
$$
 
\begin{theorem}\label{thm.llt}
For all $\kappa_1\in C(\mathcal F_\theta)$, $\kappa_2\in C(\mathcal F_{\mathrm i\theta})$ and $F\in C_c(\mathbb R\times W)$, the following hold.
\begin{enumerate}
\item
\begin{multline}\label{eq.llt1}
 \lim_{T\to\infty}\frac{\Psi(T)}{T^r}e^{-\psi(u)T} \sum_{\gamma\in\Gamma} \kappa_1(x_\gamma)\kappa_2(y_\gamma) F\left(t_\gamma-T,\frac{z_\gamma}{\sqrt T}\right)\\
 =\mathfrak j_u\,\nu(\kappa_1)\nu_{\mathrm i}(\kappa_2) \int_{\mathbb R\times W} e^{\psi(u)s}F(s,w)J_u(w)\,ds\,dw.
\end{multline}
\item
\begin{multline}\label{eq.llt2}
 \lim_{T\to\infty}\Psi(T)e^{-\psi(u)T} \sum_{\gamma\in\Gamma} \kappa_1(x_\gamma)\kappa_2(y_\gamma) F(t_\gamma-T,z_\gamma)\\
 =\mathfrak j_uJ_u(0)\nu(\kappa_1)\nu_{\mathrm i}(\kappa_2) \int_{\mathbb R\times W} e^{\psi(u)s}F(s,w)\,ds\,dw.
\end{multline}
\end{enumerate}
If local mixing has a polynomial error of order $d>0$ and $\Gamma$ is linearly $\theta$-regular near $u$ in the sense of Definition \ref{def.lin}, then there exists an integer $k\geq1$ such that the errors in \eqref{eq.llt1} and \eqref{eq.llt2} are
$$
 O_D\!\left( T^{-\eta}\|\kappa_1\|_{C^k}\|\kappa_2\|_{C^k}\|F\|_{C^1} \right), \qquad \eta=\min\left\{1,\frac{d}{2(q+k)+1}\right\},
$$
uniformly for smooth $\kappa_1,\kappa_2$ and $F\in C_c^1(\mathbb R\times W)$ supported in a fixed compact set $D\subset\mathbb R\times W$.
\end{theorem}

\subsection*{BMS correlations and partial Cartan projections}
In the Hopf parametrization \eqref{eq.hopf0}, we have
\begin{align}
 &d\widetilde m(x,y,b)=e^{\psi(2\mathcal G(x,y))}\,d\nu(x)d\nu_{\mathrm i}(y)\,db,\label{eq.bms}\\
 &g(x,y,b)=(gx,gy,b+\beta_{gx}^\theta(o,go)),\qquad (x,y,b)\exp v=(x,y,b+v).\notag
\end{align}
Fix $0<\epsilon_0\leq1$ and a nonnegative $\varphi\in C_c^\infty(\mathfrak a_\theta)$ with integral one, supported in the unit ball $B_{\mathfrak a_\te}(0,1)$ and positive near zero. 
For $0<\epsilon\leq\epsilon_0$, set
\begin{equation}\label{eq.phie}
 \varphi_\epsilon(b)=\epsilon^{-q}\varphi(b/\epsilon),\qquad k_\epsilon(h)=\int_{\mathfrak a_\theta}\varphi_\epsilon(b-h)\varphi_\epsilon(b)\,db.
\end{equation}
Then
\begin{equation}\label{eq.ke}
 \begin{gathered} k_\epsilon\geq0,\qquad \int_{\mathfrak a_\te} k_\epsilon(h)\,dh=1,\qquad \operatorname{supp}k_\epsilon\subset B_{\mathfrak a_\te}(0,2\epsilon),\\
 \|k_\epsilon\|_\infty\ll\epsilon^{-q},\qquad \norm{k_\epsilon}_{C^1} \ll\epsilon^{-q-1}. \end{gathered}
\end{equation}
 We also have $k_\epsilon(-h)=k_\epsilon(h)$. 
Fix $\kappa_1\in C(\mathcal F_\theta)$ and $\kappa_2\in C(\mathcal F_{\mathrm i\theta})$. 
Define the \emph{smoothed counting sum}
\begin{equation}\label{eq.sms}
 \mathsf S_{\epsilon}(v)=\sum_{\gamma\in\Gamma}e^{-\psi(\mu_\theta(\gamma))}\kappa_1(x_\gamma)\kappa_2(y_\gamma)k_\epsilon(v-\mu_\theta(\gamma)), \qquad (v\in \fa_\te).
\end{equation}
Since $\mu_\theta: \Gamma\to\mathfrak a_\te$ is proper, the sum is locally finite and defines a smooth function. 
When we need to indicate the dependence on $\kappa_1,\kappa_2$, we write $\mathsf S_{\epsilon}=\mathsf
S_{\epsilon}^{\kappa_1,\kappa_2}$; note that $|\mathsf S^{\kappa_1,\kappa_2}_\epsilon|\leq\|\kappa_1\|_\infty\|\kappa_2\|_\infty\mathsf S^{1,1}_\epsilon$.
 
\subsection*{Definition of $\rho$, $\rho_1$, $\rho_2$}
We fix a nonnegative $\rho\in C_c^\infty(\mathcal F_\theta^{(2)})$ such that
$$
\int\rho(x,y)\,d\nu_{\mathrm i}(y)>0 \text{ for all }x\in\Lambda_\theta,\quad \int\rho(x,y)\,d\nu(x)>0\text{ for all }y\in\Lambda_{\mathrm i\theta}.
$$
 
For each $x\in\Lambda_\theta$, choose $y\in\Lambda_{\mathrm i\theta}$ in general position with $x$, and a nonnegative $\rho_x\in C_c^\infty(\mathcal F_\theta^{(2)})$ which is positive at $(x,y)$ and supported in a product neighborhood of $(x,y)$.
Since $y\in\operatorname{supp}\nu_{\mathrm i}$, we have $\int\rho_x(x',y')\,d\nu_{\mathrm i}(y')>0$ for all $x'$ in a neighborhood of $x$. 
By compactness of $\Lambda_\theta$, a finite sum of such functions has positive $\nu_{\mathrm i}$-integral for every $x\in\Lambda_\theta$. Adding a second finite sum,
constructed with the roles of $x$ and $y$ reversed, gives a function $\rho$ satisfying both inequalities.
By continuity, there are open neighborhoods $U_1\supset\Lambda_\theta$ and $U_2\supset\Lambda_{\mathrm i\theta}$ such that
$$
 \int\rho(x,y)\,d\nu_{\mathrm i}(y)>0 \text{ for all $x\in U_1$},
 \quad
 \int\rho(x,y)\,d\nu(x)>0 \text{ for all $y\in U_2$}.
$$
Choose $\chi_1\in C_c^\infty(U_1)$ and $\chi_2\in C_c^\infty(U_2)$ with $0\leq\chi_i\leq1$, identically one on neighborhoods of $\Lambda_\theta$ and $\Lambda_{\mathrm i\theta}$, respectively.
Set
\begin{equation*}
\rho_1(x,y):=\frac{\chi_1(x)\rho(x,y)}{\int\rho(x,y')\,d\nu_{\mathrm i}(y')},\qquad
\rho_2(x,y):=\frac{\chi_2(y)\rho(x,y)}{\int\rho(x',y)\,d\nu(x')}. 
\end{equation*}
Extending by zero we obtain $\rho_1,\rho_2\in C_c^\infty(\mathcal F_\theta^{(2)})$ satisfying
\begin{equation}\label{eq.eta1}
 \int\rho_1(x,y)\,d\nu_{\mathrm i}(y)=\chi_1(x),\qquad \int\rho_2(x,y)\,d\nu(x)=\chi_2(y).
\end{equation}
 
\subsection*{Definition of $f_{i,\epsilon}$}
 
Define $f_{i,\epsilon}\in C_c(\Omega)$ by
\begin{equation}\label{eq.fie}
 f_{i,\epsilon}(\Gamma z)
 :=\sum_{\gamma\in\Gamma}\widetilde f_{i,\epsilon}(\gamma z),
 \qquad z\in\widetilde\Omega,\quad i=1,2,
\end{equation}
where \begin{align*}
 &\widetilde f_{1,\epsilon}(x,y,b)=e^{-\psi(2\mathcal G(x,y))}\kappa_1(x)\varphi_\epsilon(b)\rho_1(x,y),\\
 &\widetilde f_{2,\epsilon}(x,y,b)=e^{-\psi(2\mathcal G(x,y))}\kappa_2(y)\varphi_\epsilon(b-2\mathcal G(x,y))\rho_2(x,y).
\end{align*}
 
The  exponential factors cancel the density in \eqref{eq.bms}, so
the integrals of these  functions equal $\nu(\kappa_1)$ and $\nu_{\mathrm i}(\kappa_2)$. 
The translation by $2\mathcal G(x,y)$ in  $\widetilde f_{2,\epsilon}$ gives the expression in \eqref{eq.err1} when we unfold the correlation.

\begin{lemma}\label{lem.test}
We have
\begin{equation}\label{eq.test1}
 m(f_{1,\epsilon})=\nu(\kappa_1),\qquad m(f_{2,\epsilon})=\nu_{\mathrm i}(\kappa_2).
\end{equation}
 If $j\geq0$ is an integer and $\kappa_1,\kappa_2$ are $C^j$, then
\begin{equation}\label{eq.test2}
 \|\widetilde f_{1,\epsilon}\|_{C^j}\|\widetilde f_{2,\epsilon}\|_{C^j}\ll_j\epsilon^{-2(q+j)}\|\kappa_1\|_{C^j}\|\kappa_2\|_{C^j}.
\end{equation}
The implied constant is independent of $\kappa_1,\kappa_2,\epsilon$.
\end{lemma}
 
\begin{proof}
Unfolding the sums in \eqref{eq.fie} and using \eqref{eq.bms}, we obtain
\begin{align*}
&m(f_{1,\epsilon})=\int\kappa_1(x)\rho_1(x,y)\,d\nu(x)d\nu_{\mathrm i}(y)=\nu(\kappa_1\chi_1)=\nu(\kappa_1),\\
&m(f_{2,\epsilon})=\int\kappa_2(y)\rho_2(x,y)\,d\nu(x)d\nu_{\mathrm i}(y)=\nu_{\mathrm i}(\kappa_2\chi_2)=\nu_{\mathrm i}(\kappa_2).
\end{align*}
 The functions $\widetilde f_{i,\epsilon}$ are defined on $G/M_\theta$, not only on $\widetilde\Omega$. Their supports lie in the compact set
$$
 \left\{(x,y,b)\in\mathcal F_\theta^{(2)}\times\mathfrak a_\theta:
 (x,y)\in\supp\rho,\quad
 \|b\|\leq\epsilon_0+2\max_{\supp\rho}\|\mathcal G\|\right\},
$$
independently of $\kappa_1,\kappa_2,\epsilon$. Moreover, $\mathcal G$ and its derivatives are bounded on $\supp\rho$. 
Since $\|\varphi_\epsilon\|_{C^j}\ll_j\epsilon^{-q-j}$, the product and chain rules give
$$
 \|\widetilde f_{i,\epsilon}\|_{C^j}\ll_j\epsilon^{-q-j}\|\kappa_i\|_{C^j},\qquad i=1,2.
$$
Taking the product of these two bounds proves \eqref{eq.test2}.
\end{proof}
For $v\in\mathfrak a_\theta$, define the BMS correlation
\begin{align*}
 \mathsf C_\epsilon(v)
 &:=\int_\Omega f_{1,\epsilon}(z)
 f_{2,\epsilon}(z\exp v)\,dm(z)\notag 
\\ &= \sum_{\gamma\in\Gamma}\int_{\widetilde\Omega}\tilde f_{1,\epsilon}(\gamma z\exp(-v))\tilde f_{2,\epsilon}(z)\,d\tilde m(z).
\end{align*}
Write $z=(x,\zeta,b)\in \widetilde \Omega$ in the Hopf parametrization \eqref{eq.hopf0} and set $y=\gamma\zeta$. 
Since $\nu_{\i}$ is $(\Ga,\psi\circ\i)$-conformal \eqref{eq.conf}, we have
$$
 d\nu_{\mathrm i}(\gamma^{-1}y)=e^{-\psi(\mathrm i\beta_{\gamma^{-1}y}^{\mathrm i\theta}(o,\gamma^{-1}o))}\,d\nu_{\mathrm i}(y).
$$
By \eqref{eq.gp}, we have
$$
 2\mathcal G(x,\gamma^{-1}y)+\beta_{\gamma x}^\theta(o,\gamma o)=2\mathcal G(\gamma x,y)+\mathrm i\beta_{\gamma^{-1}y}^{\mathrm i\theta}(o,\gamma^{-1}o).
$$
Applying \eqref{eq.bms} and integrating in $b$, we obtain the following formula. 
Here and below, the integrals in $x,y$ are over $\Lambda_\theta\times\Lambda_{\mathrm i\theta}$. 
The integrand is understood to be zero unless both $(\gamma x,y)$ and $(x,\gamma^{-1}y)$ belong to $\mathcal F_\theta^{(2)}$. 
This convention is compatible with the compact supports of $\rho_1,\rho_2$.
\begin{align}
& \mathsf C_\epsilon(v)=\sum_{\gamma\in\Gamma}\int\kappa_1(\gamma x)\kappa_2(\gamma^{-1}y)\rho_1(\gamma x,y)\rho_2(x,\gamma^{-1}y)\notag\\*
 &\cdot e^{-\psi(2\mathcal G(x,\gamma^{-1}y)+\beta_{\gamma x}^\theta(o,\gamma o))} k_\epsilon\bigl(v-2\mathcal G(x,\gamma^{-1}y)-\beta_{\gamma x}^\theta(o,\gamma o)\bigr)\,d\nu(x)d\nu_{\mathrm i}(y).\label{eq.unfold}
\end{align}

Fix a closed convex cone $\mathcal C_0\subset\mathfrak a_\theta$ such that
\begin{equation}\label{eq.c0}
 u\in\operatorname{int}\mathcal C_0,
 \qquad
 \mathcal C_0\setminus\{0\}\subset\operatorname{int}\mathfrak a_\theta^+,
 \qquad \psi>0\text{ on }\mathcal C_0\setminus\{0\}.
\end{equation}
Such a cone exists by Lemma \ref{lem.lm-interior} and $\psi(u)>0$.
When \eqref{eq.rct} holds near $u$, we also choose $\mathcal C_0\setminus\{0\}\subset\mathcal U$, where $\mathcal U$ is the open cone in that condition.
 
\begin{lemma}\label{lem.locb}
For every $R>0$ there exist $t_0>0$ and $C>0$ depending on $R$ such that for all $t\geq t_0$ and all $w\in W$ with $tu+\sqrt t\,w\in\mathcal C_0$,
$$
\sum_{{\gamma\in\Gamma,
\ \|\mu_\theta(\gamma)-(tu+\sqrt t\,w)\|\leq R}}
e^{-\psi(\mu_\theta(\gamma))}\leq\frac{C}{\Psi(t)}.
$$
In particular, 
$$
 \#\{\gamma\in\Gamma:
 \|\mu_\theta(\gamma)-(tu+\sqrt t\,w)\|\leq R\}
 \ll_R \frac{e^{t\psi(u)}}{\Psi(t)}.
$$
\end{lemma}
\begin{proof}
Fix a nonnegative $\varphi_R\in C_c^\infty(\mathfrak a_\theta)$ with integral one, positive on $B_{\mathfrak a_\theta}(0,R+2)$, and define $k_R$ as in \eqref{eq.phie} using $\varphi_R$ in place of $\varphi_\epsilon$.
Construct $f_1,f_2$ as in \eqref{eq.fie}, using $\kappa_1=\kappa_2=1$ and $\varphi_R$ in place of $\varphi_\epsilon$. For $v=tu+\sqrt t\,w\in\mathcal C_0$, Definition \ref{lmp}(2) gives
\begin{equation}\label{eq.cr}
 \mathsf C_R(v):=\int_{\Omega}f_1(z)f_2(z\exp v)\,dm(z)\leq\frac{C(f_1,f_2)}{\Psi(t)}.
\end{equation}
As in \eqref{eq.unfold},
\begin{align}
 &\mathsf C_R(v)=\sum_{\gamma\in\Gamma}\int \rho_1(\gamma x,y)\rho_2(x,\gamma^{-1}y) e^{-\psi(2\mathcal G(x,\gamma^{-1}y)+\beta_{\gamma x}^\theta(o,\gamma o))}\notag\\
 &\qquad\cdot k_R\bigl(v-2\mathcal G(x,\gamma^{-1}y)-\beta_{\gamma x}^\theta(o,\gamma o)\bigr)\,d\nu(x)d\nu_{\mathrm i}(y). \label{eq.cr2}
\end{align}
Set $Q:=\supp\rho_1\cup\supp\rho_2$, and for $\ga\in\Ga$,
$$
 E_\gamma:=\{(x,y)\in\Lambda_\theta\times\Lambda_{\mathrm i\theta}:
\rho_1(x_\gamma,y)\rho_2(x,y_\gamma)>0\}.
$$
By definition, if $(x,y)\in E_\ga$, both $(x_\gamma,y)$ and $(x,y_\gamma)$ belong to $Q$. 
Applying Lemma \ref{lem.dist1} to $\gamma$ and to $\gamma^{-1}$, with $\theta$ replaced by $\mathrm i\theta$ in the latter case, gives
$$
 d_{\mathcal F_\theta}(\gamma x,x_\gamma)
 +d_{\mathcal F_{\mathrm i\theta}}(\gamma^{-1}y,y_\gamma)
 \ll_Q e^{-\gap(\gamma)}.
$$
Thus $(\gamma x,y)$ and $(x,\gamma^{-1}y)$ are uniformly close to $(x_\gamma,y)$ and $(x,y_\gamma)$, respectively, which lie in $Q$.
Therefore there exists $R_1>0$ such that both pairs lie in a fixed compact neighborhood of $Q$ in $\mathcal F_\theta^{(2)}$ whenever $\gap(\gamma)\geq R_1$.
Since $\mathcal G,\rho_1,\rho_2$ are smooth, we obtain, for all such $\gamma\in\Ga$ and $(x,y)\in E_\ga$,
\begin{align}\label{eq.ega1}
 &\bigl\|\mu_\theta(\gamma)-\beta_{\gamma x}^\theta(o,\gamma o)  - 2\mathcal G(x,\gamma^{-1}y) \bigr\|\ll e^{-\gap(\gamma)},\notag\\
 &\rho_1(\gamma x,y)\rho_2(x,\gamma^{-1}y) =\rho_1(x_\gamma,y)\rho_2(x,y_\gamma)+O(e^{-\gap(\gamma)}).
\end{align}
By \eqref{eq.eta1},
\begin{equation}\label{eq.ega2}
\int\rho_1(x_\gamma,y)\rho_2(x,y_\gamma)\,d\nu(x)d\nu_{\mathrm i}(y) =\chi_1(x_\gamma)\chi_2(y_\gamma).    
\end{equation}
We claim that there exists a finite subset $F\subset \Gamma$ such that for every $\gamma\notin F$,
\begin{enumerate}
    \item  $\chi_1(x_\gamma)=\chi_2(y_\gamma)=1$, and
$(\nu\otimes\nu_{\mathrm i})(E_\gamma)>0$.
    \item For every $(x,y)\in E_\gamma$, $\bigl\|\mu_\theta(\gamma)-\beta_{\gamma x}^\theta(o,\gamma o)- 2\mathcal G(x,\gamma^{-1}y)\bigr\|\leq 1$.
    \item  $$
 \int_{E_\gamma}\rho_1(\gamma x,y)\rho_2(x,\gamma^{-1}y)\,
 d\nu(x)d\nu_{\mathrm i}(y)\geq\frac12.
$$
\end{enumerate}
Since $\chi_1$ and $\chi_2$ are identically one on neighborhoods of $\Lambda_\theta$ and $\Lambda_{\mathrm i\theta}$, respectively, the paragraph following \eqref{eq.cf1} implies (1) outside a finite set.
For these $\gamma$, \eqref{eq.ega2} equals one. 
Its integrand vanishes off $E_\gamma$, so $E_\gamma$ has positive product measure.
By \eqref{eq.ega1} and $\theta$-regularity, (2) holds after enlarging the finite set. 
Moreover, by \eqref{eq.ega1}, \eqref{eq.ega2} and (1),
$$
 \int_{E_\gamma}
 \rho_1(\gamma x,y)\rho_2(x,\gamma^{-1}y)\,
 d\nu(x)d\nu_{\mathrm i}(y)
 =1+O(e^{-\gap(\gamma)}).
$$
Since $\gap(\gamma)\to\infty$ as $\gamma\to\infty$ in $\Gamma$, enlarging $F$ once more gives (3).

 We may assume that $F\neq\emptyset$. Choose $t_0=t_0(R)\geq1$ such that
$$
 \psi(u)t_0-\|\psi\|_*R>\max_{\gamma\in F}\psi(\mu_\theta(\gamma)).
$$
If $t\geq t_0$ and $\|\mu_\theta(\gamma)-v\|\leq R$, then
$$
 \psi(\mu_\theta(\gamma))\geq\psi(v)-\|\psi\|_*R=\psi(u)t-\|\psi\|_*R
>\max_{\gamma\in F}\psi(\mu_\theta(\gamma)),
$$
and hence $\gamma\notin F$.
Since all summands in \eqref{eq.cr2} are nonnegative, we may restrict the sum to $\gamma$ with $\|\mu_\theta(\gamma)-v\|\leq R$ and each remaining integral
to $E_\gamma$. 
For these $\gamma$, (2) gives $ \|v-2\mathcal G(x,\gamma^{-1}y) -\beta_{\gamma x}^\theta(o,\gamma o)\|\leq R+1$ for $(x,y)\in E_\gamma$.
Set $c_R:=\inf_{\|h\|\leq R+1}k_R(h)>0$.
On $E_\gamma$, the kernel factor is therefore at least $c_R$, and $e^{-\psi(2\mathcal G(x,\gamma^{-1}y) +\beta_{\gamma x}^\theta(o,\gamma o))}\geq e^{-\|\psi\|_*}e^{-\psi(\mu_\theta(\gamma))}$.
 By (3),
$$
 \mathsf C_R(v)\geq
 \frac{c_Re^{-\|\psi\|_*}}{2}
 \sum_{\|\mu_\theta(\gamma)-v\|\leq R}
 e^{-\psi(\mu_\theta(\gamma))}.
$$
Combining this with \eqref{eq.cr} proves the first claim.
Since $\psi(w)=0$, every element in the sum satisfies $ |\psi(\mu_\theta(\gamma))-t\psi(u)|\leq R\|\psi\|_*$.
Multiplying the first bound by $e^{t\psi(u)+R\|\psi\|_*}$ therefore gives the second claim.
\end{proof}
 
We now compare the BMS correlation $\mathsf C_\epsilon$ with the smoothed counting sum $\mathsf S_\epsilon$.
For $R>0$, set
$$
 \omega(R):=\sup_{\substack{\gamma\in\Gamma\\\|\mu_\theta(\gamma)\|\geq R}}e^{-\gap(\gamma)}.
$$
By $\theta$-regularity, $\omega(R)\to0$ as $R\to\infty$: for each $A>0$, the set $\{\gamma\in\Gamma:\gap(\gamma)\leq A\}$ is finite, and its partial Cartan projections are bounded.
\begin{lemma}
Let $\kappa_1\in C^\infty(\mathcal F_\theta)$, $\kappa_2\in C^\infty(\mathcal F_{\mathrm i\theta})$.
There exist $c,C,t_0>0$, independent of $\kappa_1$ and $\kappa_2$, such that for all
$0<\epsilon\leq\epsilon_0$, $t\geq t_0$ and $w\in W$ with $v=tu+\sqrt t\,w\in\mathcal C_0$,
\begin{equation}\label{eq.sece}
 \bigl|{\mathsf C}_{\epsilon}(v)-\mathsf S_{\epsilon}(v)\bigr|\leq C\epsilon^{-q-1}\omega(ct)\Psi(t)^{-1}\|\kappa_1\|_{C^1}\|\kappa_2\|_{C^1}.
\end{equation}
If \eqref{eq.rct} holds near $u$, then, after changing $c,C,t_0$ if necessary,
\begin{equation}\label{eq.sece-exp}
 \bigl|{\mathsf C}_{\epsilon}(v)-\mathsf S_{\epsilon}(v)\bigr|\leq C\epsilon^{-q-1}e^{-ct}\Psi(t)^{-1}\|\kappa_1\|_{C^1}\|\kappa_2\|_{C^1}.
\end{equation}
\end{lemma}
\begin{proof} For each $\gamma\in\Gamma$, let $A_\gamma(x,y)$ be the
integrand in the $\gamma$-term of \eqref{eq.unfold}, and let
$B_\gamma(x,y)$ be the corresponding summand of
$\mathsf S_\epsilon(v)$ multiplied by
$\rho_1(x_\gamma,y)\rho_2(x,y_\gamma)$:
\begin{align*}
A_\ga:=&\kappa_1(\gamma x)\kappa_2(\gamma^{-1}y)\rho_1(\gamma x,y)\rho_2(x,\gamma^{-1}y)\\
 &\times e^{-\psi(2\mathcal G(x,\gamma^{-1}y)+\beta_{\gamma x}^\theta(o,\gamma o))} k_\epsilon\bigl(v-2\mathcal G(x,\gamma^{-1}y)-\beta_{\gamma x}^\theta(o,\gamma o)\bigr);  \end{align*}
\begin{equation*}
B_\ga:=\kappa_1(x_\gamma)\kappa_2(y_\gamma)\rho_1(x_\gamma,y)\rho_2(x,y_\gamma)e^{-\psi(\mu_\theta(\gamma))}k_\epsilon(v-\mu_\theta(\gamma)).
\end{equation*}
It suffices to prove, uniformly for all sufficiently large $t$,
\begin{equation}\label{ab}
 |A_\gamma-B_\gamma|
 \leq C\epsilon^{-q-1}\omega(ct)e^{-\psi(\mu_\theta(\gamma))}
 \|\kappa_1\|_{C^1}\|\kappa_2\|_{C^1}.
\end{equation}
 We show below that
 \eqref{eq.eps0} holds and that $\sum_{\gamma\in\Gamma}\int B_\gamma\,d\nu\,d\nu_{\mathrm i} =\mathsf S_\epsilon(v)$ for all sufficiently large $t$.
Integrating \eqref{ab} against $\nu\otimes\nu_{\mathrm i}$ and summing over $\gamma$, we obtain
$$
 |\mathsf C_\epsilon(v)-\mathsf S_\epsilon(v)|
 \ll\epsilon^{-q-1}\omega(ct)
 \|\kappa_1\|_{C^1}\|\kappa_2\|_{C^1}
 \sum_{\|\mu_\theta(\gamma)-v\|\leq2\epsilon_0+1}
 e^{-\psi(\mu_\theta(\gamma))}.
$$
Applying Lemma \ref{lem.locb} with $R=2\epsilon_0+1$ then proves \eqref{eq.sece}. We prove \eqref{eq.sece-exp} at the end.

To prove \eqref{ab}, we first show that for all $x\in\Lambda_\theta$, $y\in\Lambda_{\mathrm i\theta}$ with $A_\gamma\ne0$ or $B_\gamma\ne0$,
\begin{equation}\label{eq.dist3}
 \begin{aligned}
 &\bigl\|2\mathcal G(x,\gamma^{-1}y)
 +\beta_{\gamma x}^\theta(o,\gamma o)
 -\mu_\theta(\gamma)\bigr\|\ll\omega(ct),\\
 &d_{\mathcal F_\theta}(\gamma x,x_\gamma)
 +d_{\mathcal F_{\mathrm i\theta}}(\gamma^{-1}y,y_\gamma)
 \ll\omega(ct).
 \end{aligned}
\end{equation}
 
First, if $A_\ga\neq 0$, then $(x,\gamma^{-1}y)\in\supp\rho_2\cap\Lambda_\theta^{(2)}$ and $ (\gamma x,y)\in\supp\rho_1\cap\Lambda_\theta^{(2)}.$
Since $\supp k_\epsilon\subset B(0,2\epsilon)$ and $\psi(v)=\psi(u) t$, we have
$$
 \left|\psi\bigl(2\mathcal G(x,\gamma^{-1}y)+\beta_{\gamma x}^\theta(o,\gamma o)\bigr)-\psi(u)t\right|\leq 2\epsilon_0\|\psi\|_{\ast}.
$$
Applying Lemma \ref{lem.dist2} with $R=\psi(u)t/2$ bounds both expressions in \eqref{eq.dist3} by $Ce^{-\gap(\gamma)}$ and gives $\|\mu_\theta(\gamma)\|\geq ct$ for all $t\geq t_1$, with $t_1$ independent of $w$ and $\epsilon$.
Thus $e^{-\gap(\gamma)}\leq\omega(ct)$, which proves \eqref{eq.dist3} in this case uniformly in $w$ and $0<\epsilon\leq\epsilon_0$.
Next, if $B_\ga\neq 0$, then $(x,y_\gamma)\in\supp\rho_2$,  $(x_\gamma,y)\in\supp\rho_1$ and $\|\mu_\theta(\gamma)-v\|\leq2\epsilon.$
Hence
$$
 \|\mu_\theta(\gamma)\|\geq\frac{\psi(u)t-2\epsilon_0\|\psi\|_*}{\|\psi\|_*}\geq ct
$$
for some fixed $c>0$ and all $t\geq t_1$ after increasing $t_1$ if necessary. Hence $e^{-\gap(\gamma)}\leq\omega(ct)$. 
Lemma \ref{lem.dist1} now gives the second bound in \eqref{eq.dist3} and shows that $(x,\gamma^{-1}y)$ lies in a fixed compact neighborhood of $\supp\rho_2$ in $\mathcal F_\theta^{(2)}$. 
Since $\mathcal G$ is smooth, we therefore have
$$
 \|\mathcal G(x,\gamma^{-1}y)-\mathcal G(x,y_\gamma)\|\ll d_{\mathcal F_{\mathrm i\theta}}(\gamma^{-1}y,y_\gamma)\ll e^{-\gap(\gamma)}.
$$
Combining this with Lemma \ref{lem.dist1} proves the remaining bound in \eqref{eq.dist3}, and hence the claim.
 
Increase $t_1$ so that the left-hand side of the first line in \eqref{eq.dist3} is at most one for all $t\geq t_1$. 
By \eqref{eq.dist3} and $\supp k_\epsilon\subset B(0,2\epsilon)$, we have
\begin{equation}\label{eq.eps0}
 \|\mu_\theta(\gamma)-v\|\leq2\epsilon_0+1
\end{equation}
whenever $t\geq t_1$ and either $A_\ga$  or $B_\ga$ is nonzero.
If $k_\epsilon(v-\mu_\theta(\gamma))\ne0$, then $\|\mu_\theta(\gamma)\|\geq ct$ as above, and hence $\chi_1(x_\gamma)=\chi_2(y_\gamma)=1$ for all $t\geq t_1$ after once more increasing $t_1$. For such $\gamma$, we therefore have
$$
 \int\rho_1(x_\gamma,y)\rho_2(x,y_\gamma)\,d\nu(x)d\nu_{\mathrm i}(y)=\chi_1(x_\gamma)\chi_2(y_\gamma)=1.
$$ 
Thus, integrating $B_\ga$ against $d\nu(x)\,d\nu_{\mathrm i}(y)$ and summing over $\gamma$ gives ${\mathsf S}_\epsilon(v)$.
The bounds \eqref{eq.dist3}, \eqref{eq.ke} and the mean value theorem give \eqref{ab}.
 
Now suppose that \eqref{eq.rct} holds near $u$. Whenever $A_\gamma$ or $B_\gamma$ is nonzero, \eqref{eq.eps0} gives $\|\mu_\theta(\gamma)-v\|\leq2\epsilon_0+1$, and the preceding argument gives $\|\mu_\theta(\gamma)\|\geq ct$.
Since $\mathcal C_0\setminus\{0\}\subset\mathcal U$ and $\mathcal C_0$ is closed, these facts imply $\mu_\theta(\gamma)\in\mathcal U$ for all sufficiently large $t$, uniformly in $w$ and $\epsilon$.
Thus \eqref{eq.rct} gives $e^{-\gap(\gamma)}\ll e^{-ct}$ after decreasing $c>0$.
The  bounds in Lemmas \ref{lem.dist1} and \ref{lem.dist2}  give \eqref{eq.dist3} and \eqref{ab} with $\omega(ct)$ replaced by $e^{-ct}$ on all nonzero summands.
The same summation and Lemma \ref{lem.locb} prove \eqref{eq.sece-exp}.
\end{proof}
 
  The
  smoothed sums satisfy the following.
\begin{proposition}\label{prop.sm}
Let $\kappa_1\in C^\infty(\mathcal F_\theta)$ and $\kappa_2\in C^\infty(\mathcal F_{\mathrm
i\theta})$.
\begin{enumerate}
\item For every fixed $0<\epsilon\leq\epsilon_0$ and $w\in W$,
\begin{equation}\label{eq.sm1}
 \lim_{t\to \infty} \Psi(t)\mathsf S_{\epsilon}(tu+\sqrt t\,w)= J_u(w)\nu(\kappa_1)\nu_{\mathrm i}(\kappa_2).
\end{equation}
\item For every fixed $0<\epsilon\leq\epsilon_0$ and $z\in W$,
\begin{equation}\label{eq.sm2}
 \lim_{t\to \infty} \Psi(t)\mathsf S_{\epsilon}(tu+z)= J_u(0)\nu(\kappa_1)\nu_{\mathrm i}(\kappa_2).
\end{equation}
\item There exist $t_0>0$ and $C>0$, independent of $\epsilon$, $\kappa_1$ and $\kappa_2$, such that for all $0<\epsilon\leq\epsilon_0$, $t\geq t_0$ and $w\in W$ with $tu+\sqrt t\,w\in\mathcal C_0$,
\begin{equation}\label{eq.sm3}
 |\Psi(t)\mathsf S_{\epsilon}(tu+\sqrt t\,w)|\leq C\epsilon^{-q}\|\kappa_1\|_\infty\|\kappa_2\|_\infty .
\end{equation}
\item Suppose in addition that \eqref{eq.rct} holds near $u$ and that the local mixing property has a polynomial error of order $d>0$. Let $k\geq1$ be the integer of Definition \ref{def.em} for a compact set $L\subset G/M_\theta$ containing the supports of all $\widetilde f_{i,\epsilon}$ (Lemma \ref{lem.test}). There exist $t_0,c,C>0$, independent of $\epsilon$, $\kappa_1$ and $\kappa_2$, such that for all $0<\epsilon\leq\epsilon_0$, $t\geq t_0$ and $w\in W$ with $tu+\sqrt t\,w\in\mathcal C_0$,
\begin{multline}\label{eq.sm4}
\bigl|\Psi(t)\mathsf S_\epsilon(tu+\sqrt t\,w)-J_u(w)\nu(\kappa_1)\nu_{\mathrm i}(\kappa_2)\bigr|\\*
 \leq C\left[\epsilon^{-2(q+k)}t^{-d}+\epsilon^{-q-1}e^{-ct}\right]\|\kappa_1\|_{C^k}\|\kappa_2\|_{C^k}.
\end{multline}
\end{enumerate}
\end{proposition}
\begin{proof}
Let $f_{1,\epsilon},f_{2,\epsilon}$ be as in \eqref{eq.fie}.
Local mixing and \eqref{eq.test1} give
$$
 \Psi(t){\mathsf C}_\epsilon(tu+\sqrt t\,w)\to J_u(w)\nu(\kappa_1)\nu_{\mathrm i}(\kappa_2).
$$
Combined with \eqref{eq.sece}, this proves (1).
Apply local mixing along $\exp(tu)$ to  the fixed function $x\mapsto
f_{2,\epsilon}(x\exp z)$  and to
$f_{1,\epsilon}$.
The first function has integral $m(f_{2,\epsilon})$ by $A_\theta$-invariance. This
gives
$$
 \Psi(t){\mathsf C}_\epsilon(tu+z)\to J_u(0)\nu(\kappa_1)\nu_{\mathrm i}(\kappa_2).
$$
Since $tu+z\in\mathcal C_0$ for sufficiently large $t$, \eqref{eq.sece} with $w=z/\sqrt t$
proves (2).
 
By \eqref{eq.ke} and Lemma \ref{lem.locb} with $R=2\epsilon_0$,
\begin{align*}
 |\mathsf S_\epsilon(v)|&\ll\epsilon^{-q}\|\kappa_1\|_\infty\|\kappa_2\|_\infty\sum_{\|\mu_\theta(\gamma)-v\|\leq2\epsilon_0}e^{-\psi(\mu_\theta(\gamma))} \\ &\ll\epsilon^{-q}\Psi(t)^{-1}\|\kappa_1\|_\infty\|\kappa_2\|_\infty,
\end{align*}
uniformly for $v=tu+\sqrt t\,w\in\mathcal C_0$ and sufficiently large $t$, proving (3).
 The supports of the smooth functions $\widetilde f_{i,\epsilon}$ lie in a fixed compact subset of $G/M_\theta$, independent of $\epsilon,\kappa_1,\kappa_2$.
Hence \eqref{eq.kd2}, \eqref{eq.test1} and \eqref{eq.test2} give
$$
 \bigl|\Psi(t){\mathsf C}_\epsilon(tu+\sqrt t\,w)-J_u(w)\nu(\kappa_1)\nu_{\mathrm i}(\kappa_2)\bigr| \ll\epsilon^{-2(q+k)}t^{-d}\|\kappa_1\|_{C^k}\|\kappa_2\|_{C^k}.
$$
Combining this with \eqref{eq.sece-exp} proves (4), since $k\geq1$.
\end{proof}

\subsection*{Proof of Theorem \ref{thm.llt}}
We first prove (1). Suppose that $\kappa_1,\kappa_2$ are smooth
and $F\in C_c^1(\mathbb R\times W)$.
Fix $R>0$ such that
$\operatorname{supp}F\subset\{(s,w):|s|\leq R,\ \|w\|\leq R\}$.
For $v=tu+z$, $z\in W$, set
$$
 F_T(v)=F\left(t-T,\frac{z}{\sqrt T}\right).
$$
For $0<\epsilon\leq\epsilon_0$, define
\begin{align*}
 L_{T,\epsilon}
 &:=\frac{\Psi(T)e^{-\psi(u)T}}{T^r}
 \sum_{\gamma\in\Gamma}\kappa_1(x_\gamma)\kappa_2(y_\gamma)
 \int_{\mathfrak a_\theta}k_\epsilon(h)e^{\psi(h)}
 F_T(\mu_\theta(\gamma)+h)\,dh\\
 &=\frac{\Psi(T)e^{-\psi(u)T}}{T^r}
 \int_{\mathfrak a_\theta}e^{\psi(v)}F_T(v)
 \mathsf S_\epsilon(v)\,dv.
\end{align*}
The second equality follows from \eqref{eq.sms}
and the change of variables $v=\mu_\theta(\gamma)+h$.

We claim that
\begin{multline}\label{eq.lte1}
 \left|L_{T,\epsilon}
 -\frac{\Psi(T)e^{-\psi(u)T}}{T^r}
 \sum_{\gamma\in\Gamma}\kappa_1(x_\gamma)\kappa_2(y_\gamma)
 F\left(t_\gamma-T,\frac{z_\gamma}{\sqrt T}\right)\right|\\*
 \ll_R\epsilon\|\kappa_1\|_\infty
 \|\kappa_2\|_\infty\|F\|_{C^1}
\end{multline}
for all sufficiently large $T$, with an implied constant independent of $\epsilon$.
 
 To prove this bound, for $C>0$ set
\begin{equation*}
 \Gamma_{C,T}:=
 \{\gamma\in\Gamma:
 |t_\gamma-T|\leq R+C\epsilon_0,\quad
 \|z_\gamma\|\leq R\sqrt T+C\epsilon_0\}.
\end{equation*}
Fix $C$ large enough that $\gamma\in\Gamma_{C,T}$ whenever $\|h\|\leq2\epsilon_0$ and either $F_T(\mu_\theta(\gamma)+h)$ or $F_T(\mu_\theta(\gamma))$ is nonzero.
Since $\dim W=2r$, the set
$$
 \left\{tu+z\in\mathfrak a_\theta:
 t\in\mathbb R, |t-T|\leq R+C\epsilon_0,
 z\in W, \|z\|\leq R\sqrt T+C\epsilon_0
 \right\}
$$
may be covered by $O_R(T^r)$ balls of radius one whose centers are of the form $t'u+z'$ with $|t'-T|\ll_R1$.
For all sufficiently large $T$, these centers lie in $\mathcal C_0$. 
Applying Lemma \ref{lem.locb} with $1$ in place of $R$ and Lemma \ref{lem.ratio}(1) gives
$$
 \sum_{\gamma\in\Gamma_{C,T}}
 e^{-\psi(\mu_\theta(\gamma))}
 \ll_R\frac{T^r}{\Psi(T)}.
$$
Since $\psi(\mu_\theta(\gamma))=\psi(u)T+O_R(1)$ for all $\gamma\in\Gamma_{C,T}$,
\begin{equation}\label{eq.gact}
 \frac{\Psi(T)e^{-\psi(u)T}}{T^r}
 \#\Gamma_{C,T}\ll_R1.
\end{equation}
Since $\|F_T\|_{C^1}\ll\|F\|_{C^1}$ for $T\geq1$ and $e^{\psi(h)}=1+O(\epsilon)$ for $\|h\|\leq2\epsilon$, the mean value theorem gives
$$
 \left|e^{\psi(h)}F_T(\mu_\theta(\gamma)+h)
 -F_T(\mu_\theta(\gamma))\right|
 \ll\epsilon\|F\|_{C^1}.
$$
Integrating against $k_\epsilon$ and summing over $\Gamma_{C,T}$ proves \eqref{eq.lte1} by \eqref{eq.gact}.
 
By \eqref{eq.lte1}, it remains to show that, for each fixed $0<\epsilon\leq\epsilon_0$,
\begin{equation*}
 \lim_{T\to\infty}L_{T,\epsilon}
 =\mathfrak j_u\,\nu(\kappa_1)\nu_{\mathrm i}(\kappa_2)
 \int_{\mathbb R\times W}e^{\psi(u)s}F(s,w)J_u(w)\,ds\,dw.
\end{equation*}
 
For all sufficiently large $T$, $\supp F_T\subset\mathcal C_0$.
Let $s=t-T$ and $w=z/\sqrt{T+s}$, so that $s\in[-R,R]$ whenever $F_T(v)\ne0$, $v=tu+z=(T+s)u+\sqrt{T+s}\, w$, and \eqref{eq.sm1} applies with $w$ fixed.
For $T>2R$, \eqref{eq.qr} and \eqref{eq.ju} give
\begin{equation}\label{eq.lti1}
 \begin{aligned}
 &L_{T,\epsilon}
 =\mathfrak j_u\int_{[-R,R]\times W}
 e^{s\psi(u)}F\left(s,\sqrt{1+\frac{s}{T}}\,w\right)
 \frac{\Psi(T)}{\Psi(T+s)}\\
 &\qquad\cdot\left(1+\frac{s}{T}\right)^r
 \bigl[\Psi(T+s)\mathsf S_\epsilon
 ((T+s)u+\sqrt{T+s}\,w)\bigr]\,ds\,dw.
 \end{aligned}
\end{equation}
For fixed $\epsilon$, these integrands are uniformly bounded and supported in one fixed compact set, by \eqref{eq.sm3} and Lemma \ref{lem.ratio}(1).
Applying \eqref{eq.sm1} at each fixed $(s,w)$ and then dominated convergence gives the claimed limit.
Letting $T\to\infty$ and then $\epsilon\to0$ in \eqref{eq.lte1} proves \eqref{eq.llt1} for these functions.
 
Approximating $\kappa_1,\kappa_2,F$ uniformly by smooth functions, with the supports of the smooth functions approximating $F$ in a fixed compact set, proves (1) for continuous functions, using \eqref{eq.gact} and boundedness of $J_u$.
For the quantitative statement, assume that local mixing has a polynomial error of order $d$
and that $\Gamma$ is linearly $\theta$-regular near $u$ in the sense of \eqref{eq.rct}, and let $k$ be as in \eqref{eq.sm4}.
On the fixed compact set of integration in \eqref{eq.lti1}, \eqref{eq.sm4}, Lemma \ref{lem.ratio}(3) and $|(1+s/T)^r-1|\ll T^{-1}$ give
\begin{multline*}
 \left|L_{T,\epsilon}
 -\mathfrak j_u\nu(\kappa_1)\nu_{\mathrm i}(\kappa_2)
 \int e^{\psi(u)s}F(s,w)J_u(w)\,ds\,dw\right|\\*
 \ll_D\left[T^{-1}+T^{-d}
 +\epsilon^{-2(q+k)}T^{-d}
 +\epsilon^{-q-1}e^{-cT}\right]
 \|\kappa_1\|_{C^k}\|\kappa_2\|_{C^k}\|F\|_{C^1}.
\end{multline*}
 Here replacing $F(s,\sqrt{1+s/T}\,w)$ by $F(s,w)$ changes the integral by $O_D(T^{-1}\|F\|_{C^1})$. 
 Combining this bound with \eqref{eq.lte1} and choosing $\epsilon=T^{-d/(2(q+k)+1)}$ proves the claimed bounds on the error term.
 
We now prove (2). 
First suppose that $\kappa_1,\kappa_2$ are smooth and $F\in C_c^1(\mathbb R\times W)$, and use the function $\mathsf S_\epsilon$ defined by
\eqref{eq.sms}.
Repeating the preceding argument with \eqref{eq.sm2} in place of \eqref{eq.sm1}, we obtain, for fixed
$0<\epsilon\leq\epsilon_0$,
\begin{multline*}
 \lim_{T\to\infty}\Psi(T)e^{-\psi(u)T}
 \int_{\mathfrak a_\theta}
 e^{\psi(v)}F(t-T,z)\mathsf S_\epsilon(v)\,dv\\
 =\mathfrak j_uJ_u(0)\nu(\kappa_1)\nu_{\mathrm i}(\kappa_2)
 \int_{\mathbb R\times W}e^{\psi(u)s}F(s,z)\,ds\,dz,
 \qquad v=tu+z.
\end{multline*}
Here we use $t=T+s$, $dv=\mathfrak j_u\,dt\,dz$ and dominated convergence on a fixed compact set.
The definition of $\mathsf S_\epsilon$, the bounds $\|h\|\leq2\epsilon$ on $\supp k_\epsilon$ and $\int k_\epsilon=1$, and Lemma \ref{lem.locb} show that replacing the normalized integral by
$$
 \Psi(T)e^{-\psi(u)T}
 \sum_{\gamma\in\Gamma}\kappa_1(x_\gamma)\kappa_2(y_\gamma)
 F(t_\gamma-T,z_\gamma)
$$
changes it by $O_F(\epsilon\|\kappa_1\|_\infty\|\kappa_2\|_\infty)$.
Letting $T\to\infty$ and then $\epsilon\to0$ proves \eqref{eq.llt2} for these functions.
Uniform approximation, with errors controlled by Lemma \ref{lem.locb} as above, gives the result for continuous $\kappa_1,\kappa_2,F$. 
For the error term in (2), assume again that local mixing has a polynomial error
of order $d$ and that $\Gamma$ is linearly $\theta$-regular near $u$. We first show that, for every compact set $B\subset W$,
$$
 \sup_{z\in B}|J_u(z/\sqrt t)-J_u(0)|\ll_B t^{-d}.
$$
Fix a nonnegative $h\in C_c(\Omega)$ obtained from a smooth compactly supported lift as in \eqref{eq.per}, with $m(h)>0$.
Apply \eqref{eq.kd2} to $\int_\Omega h(x)h(x\exp(tu+z))\,dm(x)$ in two ways: with $w=z/\sqrt t$, and with $w=0$ after replacing the translated function by $x\mapsto h(x\exp z)$. Since $m(h(\,\cdot\,\exp z))=m(h)$, subtracting the two formulas and dividing by $m(h)^2>0$ gives
$J_u(z/\sqrt t)=J_u(0)+O_B(t^{-d})$, uniformly for $z\in B$.
 The translated lifts have uniformly bounded smooth norms and supports in a common compact set.
 
 Together with \eqref{eq.sm4}, this proves
\begin{multline}\label{eq.tube1}
 \left|\Psi(t)\mathsf S_\epsilon(tu+z)
 -J_u(0)\nu(\kappa_1)\nu_{\mathrm i}(\kappa_2)\right|\\*
 \ll_B\left[\epsilon^{-2(q+k)}t^{-d}
       +\epsilon^{-q-1}e^{-ct}\right]
       \|\kappa_1\|_{C^k}\|\kappa_2\|_{C^k},\qquad z\in B,
\end{multline}
for all sufficiently large $t$, uniformly in $0<\epsilon\leq\epsilon_0$.
Here $k$ is the integer in \eqref{eq.sm4}. The integer from
Definition \ref{def.em} used for the translates of $h$ only affects the implied constant, since $h$ is fixed. So $k$ does not depend on $B$.
Integrating over the fixed support of $F$ and using Lemma \ref{lem.ratio}(3) gives the bound for the smoothed integral in (2), with no $T^{-1}$ term. The analogue of
\eqref{eq.lte1} for $F(t-T,z)$, followed by $\epsilon=T^{-d/(2(q+k)+1)}$, proves the error bound in (2).

\subsection*{Counting in tubes}
Chow-Oh \cite[Theorem 6.1]{CO} proved the corresponding asymptotic for Borel Anosov groups. 
We use $\delta_u=\psi(u)$ for the direction fixed in this section.
\begin{theorem}\label{thm.tubes}
Let $D\subset W$ be compact with positive Lebesgue measure and boundary of measure zero, and let $b:D\to\mathbb R$ be continuous. 
Set
$$
 \mathbb T_{T,b}:=\{tu+z:z\in D,\ 0\leq t\leq T+b(z)\}.
$$
Then, for all $\kappa_1\in C(\mathcal F_\theta)$ and $\kappa_2\in C(\mathcal F_{\i\theta})$, as $T\to\infty$,
\begin{multline*}
 \sum_{{\gamma\in\Gamma,\, 
                  \mu_\theta(\gamma)\in\mathbb T_{T,b}}}
 \kappa_1(x_\gamma)\kappa_2(y_\gamma)\\
 =\left(
 \frac{\mathfrak j_uJ_u(0)}{\delta_u}
 \nu(\kappa_1)\nu_{\mathrm i}(\kappa_2)
 \int_D e^{\delta_u b(z)}\,dz+o(1)
 \right)\frac{e^{\delta_uT}}{\Psi(T)}.
\end{multline*}
If, in addition, local mixing has a polynomial error, $\Gamma$ is linearly $\theta$-regular near $u$ in the sense of \eqref{eq.rct}, $D$ is the closure of a bounded domain with piecewise $C^1$ boundary and $b$ is Lipschitz, then the term $o(1)$ can be replaced by $O_{D,b}(T^{-\eta'}\|\kappa_1\|_{C^k}\|\kappa_2\|_{C^k})$ for smooth $\kappa_1,\kappa_2$ and some $\eta'>0$ and integer $k\geq1$.
\end{theorem}
\begin{proof}
For every sufficiently large $n$, Lemma \ref{lem.locb} gives
$$
 \#\{\gamma\in\Ga:z_\gamma\in D,\ n\leq t_\gamma<n+1\}
 \ll_D\frac{e^{\delta_un}}{\Psi(n)}.
$$
Each Cartan projection counted on the left lies in a ball of radius $\|u\|+\max_{z\in D}\|z\|$ centered at $nu$. 
Summing over $n$ and applying Lemma \ref{lem.ratio}(2) with $\eta=\delta_u/2$ gives
$$
 \limsup_{T\to\infty}\Psi(T)e^{-\delta_uT}
 \#\{\gamma:z_\gamma\in D,\ 0\leq t_\gamma<T-R\}
 \ll_D e^{-\delta_uR/2}.
$$
The finitely many elements with bounded $t_\gamma$ and $z_\gamma\in D$ contribute zero in the limit, since $\log\Psi(T)=o(T)$.
The same sum with factors $\kappa_1(x_\gamma)\kappa_2(y_\gamma)$ is bounded in absolute value by the count above times $\|\kappa_1\|_\infty\|\kappa_2\|_\infty$.
 
For $R>\max\{0,-\min_D b\}$,  the compact set $\{(s,z)\in\mathbb R\times W:z\in D,\ -R\leq s\leq b(z)\}$ has boundary of Lebesgue measure zero. 
Its boundary is contained in the two graphs $s=-R$, $s=b(z)$ over $D$, together with a bounded subset of $\mathbb R\times\partial D$. 
Fubini's theorem gives the claim.  Applying \Cref{thm.llt}(2) to continuous upper and lower approximations to this
indicator gives the following limit. For arbitrary
$\kappa_1,\kappa_2$, the approximation error is bounded by $\|\kappa_1\|_\infty\|\kappa_2\|_\infty$ times the  corresponding error for $\kappa_1=\kappa_2=1$:
\begin{multline*}
 \lim_{T\to\infty}\Psi(T)e^{-\delta_uT}
 \sum_{\substack{\gamma\in\Gamma,\ z_\gamma\in D\\
                  T-R\leq t_\gamma\leq T+b(z_\gamma)}}
 \kappa_1(x_\gamma)\kappa_2(y_\gamma)\\
 =\mathfrak j_uJ_u(0)\nu(\kappa_1)\nu_{\mathrm i}(\kappa_2)
 \int_D\int_{-R}^{b(z)}e^{\delta_us}\,ds\,dz.
\end{multline*}
Letting $R\to\infty$ and using the preceding tail bound proves the claimed limit.
 
For the version with error term, suppose that local mixing has a polynomial error of order $d>0$, that $\Gamma$ is linearly $\theta$-regular near $u$, that $D$ has piecewise $C^1$ boundary, and that $b$ is Lipschitz.
Extend $b$ to a Lipschitz function on a fixed neighborhood of $D$. 
First assume that $\kappa_1,\kappa_2$ are smooth and nonnegative.
 
For a bounded Borel set $E\subset\mathfrak a_\theta$,
\eqref{eq.sms} and the change of variables $v=\mu_\theta(\gamma)+h$ give
\begin{align*}
 &\int_E e^{\psi(v)}\mathsf S_\epsilon(v)\,dv\\
 &\quad=\sum_{\gamma\in\Gamma}\kappa_1(x_\gamma)\kappa_2(y_\gamma)
 \int_{\mathfrak a_\theta}k_\epsilon(h)e^{\psi(h)}
 \mathbf1_E(\mu_\theta(\gamma)+h)\,dh.
\end{align*}
This identity and $\supp k_\epsilon\subset B(0,2\epsilon)$ bound the count between smoothed integrals over slightly smaller and larger tubes, up to factors $e^{\pm C\epsilon}$.
More precisely, for a sufficiently large fixed $C$, replace $D$ by $D_\epsilon^-:=\{z\in D: \operatorname{dist}(z,W\setminus D)\geq C\epsilon\}$ and $D_\epsilon^+:=\{z\in W: \operatorname{dist}(z,D)\leq C\epsilon\}$, respectively.
The smaller tube has interval in the $t$-coordinate $[C\epsilon,T+b(z)-C\epsilon]$, and the larger tube has interval $[-C\epsilon,T+b(z)+C\epsilon]$.
 
We apply \eqref{eq.tube1} to these integrals. 
The sets $D_\epsilon^\pm$ all lie in one fixed compact subset of $W$.
Choose $t_0$ large enough for that bound and for Lemma \ref{lem.ratio}(2). 
The integral over $t<t_0$ is $O_{D,b}(\|\kappa_1\|_\infty\|\kappa_2\|_\infty)$, uniformly in $\epsilon$, by the identity above, since $\mu_\theta|_\Gamma$ is proper. 
After multiplication by $\Psi(T)e^{-\delta_uT}$,  this is
$O_{D,b}(e^{-cT}\|\kappa_1\|_\infty\|\kappa_2\|_\infty)$ for some $c>0$. On $t_0\leq t<T/2$,  \eqref{eq.sm3} and Lemma \ref{lem.ratio}(2), with $\eta=\delta_u/2$, give
$$
 \epsilon^{-q}\int_{t_0}^{T/2}
 e^{-\delta_u(T-t)}\frac{\Psi(T)}{\Psi(t)}\,dt
 \ll \epsilon^{-q}e^{-\delta_uT/4}.
$$
The normalized contribution from this interval is therefore
$$ O_{D,b}\!\left(\epsilon^{-q}e^{-cT}\|\kappa_1\|_\infty\|\kappa_2\|_\infty\right)
$$ after decreasing $c>0$.
 
For $T/2\leq t\leq T+b(z)+C\epsilon$, set $s=T-t$.
Here $s$ is bounded below independently of $T,\epsilon$, and Lemma \ref{lem.ratio}(3) gives
$$
 \left|\frac{\Psi(T)}{\Psi(T-s)}-1\right|
 \ll_{D,b} (|s|+1)T^{-d}e^{C_\Psi(|s|+1)T^{-d}}.
$$
This is integrable against $e^{-\delta_us}\,ds$ on this range, with integral $O_{D,b}(T^{-d})$. 
Using $dv=\mathfrak j_u\,dt\,dz$ and \eqref{eq.tube1}, the normalized smaller and larger integrals have main terms
$$
 \frac{\mathfrak j_uJ_u(0)}{\delta_u}
 \nu(\kappa_1)\nu_{\mathrm i}(\kappa_2)
 \int_{D_\epsilon^\pm}
      e^{\delta_u(b(z)\pm C\epsilon)}\,dz,
$$
respectively, and errors bounded by
$$
 C_{D,b}\left[\epsilon^{-2(q+k)}T^{-d}
       +\epsilon^{-q-1}e^{-cT}\right]
       \|\kappa_1\|_{C^k}\|\kappa_2\|_{C^k}.
$$
The exponentially small terms above are smaller than this bound.
 
Since $D$ has piecewise $C^1$ boundary, $\operatorname{vol}_W(D_\epsilon^+\setminus
D_\epsilon^-)=O_D(\epsilon)$. 
The two main terms therefore differ from the stated main term by $O_{D,b}(\epsilon\|\kappa_1\|_\infty\|\kappa_2\|_\infty)$.
The factors $e^{\pm C\epsilon}$ also change the main term by $O_{D,b}(\epsilon\|\kappa_1\|_\infty\|\kappa_2\|_\infty)$.
Taking $\epsilon=T^{-d/(2(q+k)+1)}$ therefore proves the conclusion, with $\eta'=d/(2(q+k)+1)$.
For real smooth $\kappa_i$, write $\kappa_i=(\kappa_i+\|\kappa_i\|_\infty)-\|\kappa_i\|_\infty$.
Each of the two nonnegative functions has $C^k$ norm at most $2\|\kappa_i\|_{C^k}$.
Expanding the product and then separating real and imaginary parts proves the general case.
 
The constants in this proof are uniform for Euclidean balls in $W$ whose centers lie in a fixed compact set and whose radii lie in a fixed compact subinterval of $(0,\infty)$, and for $b$ in a family with uniformly bounded sup norms and Lipschitz constants.
The supports used in the proof then remain in one compact set, and $\operatorname{vol}_W(D_\epsilon^+\setminus D_\epsilon^-)=O(\epsilon)$ uniformly.
When $W=\{0\}$, take $D=\{0\}$ with mass one; then no boundary term arises.
\end{proof}
 
To deduce Theorem \ref{thm.cnt1}(2), now take the unit maximizing vector $u$ from the introduction, so that $W=u^\perp$ and $\mathfrak j_u=1$. 
For any $0<\e'<\e$ and $a>0$, the count in Theorem \ref{thm.cnt1}(2) lies, for all sufficiently large $T$, between the counts in
$$
 [0,T-a]u+\{z\in W:\|z\|\leq\epsilon'\}
 \quad\text{and}\quad
 [0,T]u+\{z\in W:\|z\|\leq\epsilon\},
$$
since $\|tu+z\|^2=t^2+\|z\|^2$. Apply Theorem \ref{thm.tubes} with $\kappa_1=\kappa_2=1$, taking $b\equiv-a$ for the lower bound and $b\equiv0$ for the upper bound.
Letting $T\to\infty$, then $a\to 0$ and $\epsilon'\to \epsilon$, proves Theorem \ref{thm.cnt1}(2).

\section{Counting in Euclidean balls}\label{sec.eucl}
Let $\langle\cdot,\cdot\rangle$ be the Weyl-group invariant inner product and $\|\cdot\|$ be the induced norm fixed in Section \ref{ps}. 
Set
$$
 \delta:=\max_{\substack{v\in\mathfrak a_\theta^+, \|v\|=1}}\psi_\Gamma^\theta(v)>0.
$$
 The maximum exists by upper semicontinuity of $\psi_\Gamma^\theta$ \cite[Theorem 3.3]{KOW}, and $\delta>0$ since $\psi_\Gamma^\theta>0$ on $\operatorname{int}\mathcal L_\theta\neq\emptyset$.
\begin{lemma}\label{lem.dir}
There is a unique unit vector $u\in\mathfrak a_\theta^+$ such that $\psi_\Gamma^\theta(u)=\delta$.
The linear form $\psi(v):=\delta\langle u,v\rangle$ is tangent to $\psi_\Gamma^\theta$ at $u$ and is $(\Gamma,\theta)$-proper. 
If $\Gamma$ itself is $\theta$-Anosov, then $u\in\operatorname{int}\mathcal L_\theta$.
\end{lemma}
 \begin{proof}
Choose a unit vector $u\in \fa_\te^+$ with $\psi_\Gamma^\theta(u)=\delta$.
Concavity and homogeneity give $\delta+s\psi_\Gamma^\theta(v) \leq\psi_\Gamma^\theta(u+sv)\leq\delta\|u+sv\|$ for $s>0$ whenever $\psi_\Gamma^\theta(v)>-\infty$. 
Hence
$$
 \psi_\Gamma^\theta(v)
 \leq\delta\lim_{s\to 0}\frac{\|u+sv\|-1}{s}
 =\delta\langle u,v\rangle.
$$
The inequality is automatic where $\psi_\Gamma^\theta=-\infty$, so $\psi$ is tangent at $u$. 
If another unit vector $w$ attains the maximum, then $\delta\leq\delta\langle u,w\rangle\leq\delta$, which implies $w=u$.
Since $\psi_\Gamma^\theta(u)>0$, we have $u\in\mathcal L_\theta\setminus\{0\}$.
Let $\{\omega_\alpha^\vee:\alpha\in\Pi\}$ be the basis of $\mathfrak a$ dual to the simple roots, so that $\beta(\omega_\alpha^\vee)=\delta_{\alpha\beta}$.
The entries $\langle\omega_\alpha^\vee,\omega_\beta^\vee\rangle$ are nonnegative, since the inverse of the Cartan matrix has nonnegative entries.
Since $u\in\mathfrak a_\theta^+\setminus\{0\}$,  there is $\alpha\in\theta$ with $\alpha(u)>0$.
Expanding $u$ and $H\in\mathfrak a^+$ in this basis gives
$$
 \langle u,H\rangle\geq\alpha(u)\alpha(H)\|\omega_\alpha^\vee\|^2.
$$
Since $p_\theta$ is the orthogonal projection and $u\in\mathfrak a_\theta$, it follows that
$$
 \psi(\mu_\theta(\gamma))
 =\delta\langle u,\mu(\gamma)\rangle
 \geq\delta\alpha(u)\|\omega_\alpha^\vee\|^2\alpha(\mu(\gamma))
 \to+\infty
$$
as $\gamma\to\infty$ in $\Gamma$, by $\theta$-regularity.
Thus $\psi$ is $(\Gamma,\theta)$-proper.

If $\Gamma$ is $\theta$-Anosov,  then, since $\psi$ is a $(\Gamma,\theta)$-proper tangent form at $u$, \cite[Theorem 13.2]{KOW} (see also \cite[Corollary 5.9.1]{Samb4}) gives $u\in\operatorname{int}\mathcal
L_\theta$.
\end{proof}
 
Fix PS measures $\nu=\nu_\psi$ and $\nu_{\mathrm i}=\nu_{\psi\circ\mathrm i}$, and assume that their BMS measure satisfies Definition \ref{lmp} with respect to this $u$. Denote the functions in that definition by $\Psi$ and $J_u$.
Lemma \ref{lem.lm-interior} gives $u\in\operatorname{int}\mathcal L_\theta\subset\operatorname{int}\mathfrak a_\theta^+$.
We apply the results of Section \ref{sec.pc} with these choices of $u,\psi,\nu,\nu_{\mathrm i}$. 
In particular,
\begin{equation}\label{eq.cov}
 W=u^\perp,\qquad dv=t^r\,dt\,dw,\qquad \psi(tu+\sqrt t\,w)=\delta t,
\end{equation}
where $t>0$, $w\in W$, and $v=tu+\sqrt t\,w$.
For $\kappa_1\in C(\mathcal F_\theta)$ and $\kappa_2\in C(\mathcal F_{\mathrm i\theta})$, set
\begin{equation*}
 \mathcal N_{\kappa_1,\kappa_2}(T):=\sum_{\substack{\gamma\in\Gamma,  \|\mu_\theta(\gamma)\|<T}}\kappa_1(x_\gamma)\kappa_2(y_\gamma).
\end{equation*}

\begin{lemma}\label{lem.subexp}
For every open cone $\mathcal U\subset\mathfrak a_\theta^+$ containing $u$, there is $\zeta\in(0,\delta)$ such that
\begin{equation*}
 \#\{\gamma:\|\mu_\theta(\gamma)\|<T,\ \mu_\theta(\gamma)\notin\mathcal U\}\ll e^{(\delta-\zeta)T}.
\end{equation*}
\end{lemma}
 
\begin{proof}
By Lemma \ref{lem.dir}, $\psi_\Gamma^\theta(v)\leq\delta\langle u,v\rangle<\delta$ for all unit vectors $v \in\mathfrak a_\te^+-\cal U$.
 The unit vectors in $\mathfrak
 a_\theta^+\setminus\mathcal U$ form a compact set disjoint from $u$. Hence there is
 $a\in(0,\delta)$  such that
 $\psi_\Gamma^\theta(v)<a$ for every such vector. 
By the definition of $\psi_\Ga^\te$, each such direction has an open cone $\mathcal U_j$ around it for which $\sum_{\mu_\theta(\gamma)\in\mathcal U_j}e^{-a\|\mu_\theta(\gamma)\|}<\infty$.
A finite collection covers the compact set, and
$$
 \#\{\gamma:\mu_\theta(\gamma)\in\mathcal U_j,\ \|\mu_\theta(\gamma)\|<T\}\leq e^{aT}\sum_{\mu_\theta(\gamma)\in\mathcal U_j}e^{-a\|\mu_\theta(\gamma)\|}\ll e^{aT}.
$$
There are only finitely many $\gamma\in\Gamma$ with $\mu_\theta(\gamma)=0$, since $\mu_\theta|_\Gamma$ is proper. 
This proves the lemma with $\zeta=\delta-a$. 
If the compact set is empty, the same conclusion follows directly from properness.
\end{proof}
 
Recall $k_\epsilon$ and $\mathsf S_\epsilon$ from \eqref{eq.phie} and \eqref{eq.sms}. 
Set
\begin{equation*}
 I_{T,\epsilon}:=\int_{\|v\|<T}e^{\psi(v)}\mathsf S_\epsilon(v)\,dv.
\end{equation*}
 
\begin{lemma}\label{lem.ntit}
For any nonnegative $\kappa_1\in C(\mathcal F_\theta)$, $\kappa_2\in C(\mathcal F_{\mathrm
i\theta})$ and $T>2\epsilon>0$,
\begin{equation}\label{eq.ntit}
 e^{-2\delta\epsilon}I_{T-2\epsilon,\epsilon}\leq\mathcal N_{\kappa_1,\kappa_2}(T)\leq e^{2\delta\epsilon}I_{T+2\epsilon,\epsilon}.
\end{equation}
\end{lemma}
 
\begin{proof}
Interchanging the sum and integral and setting $h=v-\mu_\theta(\gamma)$ gives
\begin{equation}\label{eq.ite}
 I_{T,\epsilon}=\sum_{\gamma\in\Gamma}\kappa_1(x_\gamma)\kappa_2(y_\gamma)\int k_\epsilon(h)e^{\psi(h)}\mathbf1_{\{\|\mu_\theta(\gamma)+h\|<T\}}\,dh.
\end{equation}
On $\operatorname{supp}k_\epsilon$, we have $\|h\|\leq2\epsilon$ and $|\psi(h)|\leq2\delta\epsilon$.
The triangle inequality and $\int k_\epsilon=1$ prove the claim.
\end{proof}
 
Let $\mathcal C_0$ be as in \eqref{eq.c0}.
For the quantitative statements, we also assume that $\Gamma$ is linearly $\theta$-regular near $u$ in the sense of \eqref{eq.rct}, and choose $\mathcal C_0$ so that $\mathcal C_0\setminus\{0\}$ is contained in the open cone in that condition.
\begin{lemma}\label{lem.subexp2}
There exists $\zeta\in(0,\delta)$ such that, for each fixed $t_0>0$, as $T\to\infty$,
\begin{align}\label{eq.ite2}
 I_{T,\epsilon}
 &=\int_{t_0}^T e^{\delta t}t^r
 \int_W\mathbf1_{\mathcal C_0\cap B(0,T)}
 (tu+\sqrt t\,w)\mathsf S_\epsilon(tu+\sqrt t\,w)\,dw\,dt
 \notag\\
 &\quad+O\!\left(
 e^{(\delta-\zeta)T}
 \|\kappa_1\|_\infty\|\kappa_2\|_\infty
 \right).
\end{align}
The exponent $\zeta$ is independent of $t_0$. 
The implied constant may depend on $t_0$, $\mathcal C_0$ and $\epsilon_0$, but not on $T$, $0<\epsilon\leq\epsilon_0$ or $\kappa_1,\kappa_2$.
\end{lemma}
 
\begin{proof}
Choose an open convex cone $\mathcal C$ containing $u$ such that $\overline{\mathcal C}\setminus\{0\} \subset\operatorname{int}\mathcal C_0$.
There is $R_0>0$ such that $v\in\mathcal C$, $\|v\|\geq R_0$ and $\|h\|\leq2\epsilon_0$ imply $v+h\in\mathcal C_0$.
 
By \eqref{eq.ite}, a term contributing to the integral outside $\mathcal C_0$ must satisfy $\|\mu_\theta(\gamma)\|<T+2\epsilon_0$ and either $\|\mu_\theta(\gamma)\|<R_0$ or
$\mu_\theta(\gamma)\notin\mathcal C$.
Each term contributes at most $e^{2\delta\epsilon_0} \|\kappa_1\|_\infty\|\kappa_2\|_\infty$ in absolute value. 
Properness and Lemma \ref{lem.subexp}, applied at radius $T+2\epsilon_0$,
bound their sum by the stated error.
 
The set $\{v\in\mathcal C_0:\langle v,u\rangle<t_0\}$ is bounded, and its contribution is $O_{t_0}(\|\kappa_1\|_\infty\|\kappa_2\|_\infty)$ by \eqref{eq.ite}, since $\mu_\theta|_\Gamma$ is proper.
Applying \eqref{eq.cov} to the remaining integral proves the lemma.
\end{proof}

 \begin{lemma}\label{lem.main}
Let $\kappa_1,\kappa_2$ be smooth.
\begin{enumerate}
\item For every fixed $0<\epsilon\leq\epsilon_0$,
$$
 \lim_{T\to\infty}
 \frac{\Psi(T)e^{-\delta T}}{T^r}
 \int_{\|v\|<T}e^{\psi(v)}\mathsf S_\epsilon(v)\,dv
 =c_u\nu(\kappa_1)\nu_{\mathrm i}(\kappa_2).
$$
 
\item Suppose that $\Gamma$ is linearly $\theta$-regular near $u$ in the sense of \eqref{eq.rct} and that local mixing has a polynomial error  of order $d>0$. Let $k$ be as in Proposition \ref{prop.sm}(4).
There exists $c>0$ such that, for all sufficiently large $T$, uniformly in $0<\epsilon\leq\epsilon_0$,
\begin{multline*}
 \left|
 \frac{\Psi(T)e^{-\delta T}}{T^r}
 \int_{\|v\|<T}e^{\psi(v)}\mathsf S_\epsilon(v)\,dv
 -c_u\nu(\kappa_1)\nu_{\mathrm i}(\kappa_2)
 \right|\\*
 \ll
 \left[T^{-1}+\epsilon^{-2(q+k)}T^{-d}
 +\epsilon^{-q-1}e^{-cT}\right]
 \|\kappa_1\|_{C^k}\|\kappa_2\|_{C^k}.
\end{multline*}
\end{enumerate}
\end{lemma}
 
 \begin{proof} 
 Set
\begin{equation*}
 \lambda_T:={\Psi(T)}{T^{-r}}e^{-\delta T}.
\end{equation*}
Fix $0<\eta<\delta$ and let $t_\eta$ be as in Lemma \ref{lem.ratio}(2).
Choose $t_0>t_\eta$ sufficiently large so that \eqref{eq.sm3} in Proposition \ref{prop.sm}(3) holds for all $t>t_0$. 
Since $\log\Psi(T)=o(T)$ by Lemma \ref{lem.ratio}(2), for every $\zeta>0$ and all large $T$,
$\lambda_Te^{(\delta-\zeta)T}\leq e^{-\zeta T/2}$.
Since $\|tu+\sqrt t\,w\|^2=t^2+t\|w\|^2$, this together with \eqref{eq.ite2} gives
\begin{align}\label{eq.main1}
\lambda_TI_{T,\epsilon}
&=\int_{t_0}^T\int_W\mathbf1_{D_T}(t,w)e^{-\delta(T-t)}
\frac{\Psi(T)}{\Psi(t)}\left(\frac tT\right)^r\notag\\*
&\cdot[\Psi(t)\mathsf S_\epsilon(tu+\sqrt t\,w)]\,dw\,dt
+O\!\left(e^{-\zeta T/2}\|\kappa_1\|_\infty\|\kappa_2\|_\infty\right),
\end{align}
where $D_T=\{(t,w)\in(0,\infty)\times W: tu+\sqrt t\,w\in\mathcal C_0,\ t^2+t\|w\|^2<T^2\}.$
Since $\langle u,v\rangle\gg\|v\|$ for $v\in\mathcal C_0$, the cone condition implies $\|w\|\ll\sqrt t$.
Thus, since $\dim W=2r$, for each $t$, $\{w\in W:(t,w)\in D_T\}$ has volume $O(t^r)$.
Together with \eqref{eq.sm3} and Lemma \ref{lem.ratio}(2), this shows that the integral in \eqref{eq.main1} over $t_0\leq t<T/2$ is
\begin{equation}\label{eq.main2}
 O\!\left(\epsilon^{-q}e^{-cT}
 \|\kappa_1\|_\infty\|\kappa_2\|_\infty\right)
\end{equation}
for some $c>0$.
 
For $T/2\leq t<T$, substitute $s=T-t$ and let
\begin{align*}
& R_T:=\left\{(s,w)\in(0,\tfrac{T}2]\times W:
 \begin{gathered}
 (T-s)u+\sqrt{T-s}\,w\in\mathcal C_0,
 \|w\|^2<2s+\tfrac{s^2}{T-s}
 \end{gathered}\right\}\\
&M_T(s):=\Psi(T)\Psi(T-s)^{-1}(1-\tfrac{s}T)^r.
\end{align*}
The integral in \eqref{eq.main1} over $T/2\leq t<T$ can therefore be written as 
\begin{equation}\label{eq.main3}
\begin{split}
 I'_{T,\epsilon}
 &:=\int_0^{\frac T2}\int_W\mathbf1_{R_T}(s,w)e^{-\delta s}M_T(s) [\Psi(T-s)\mathsf S_\epsilon
 ((T-s)u+\sqrt{T-s}\,w)]\,dw\,ds.
\end{split}
\end{equation}
For fixed $\epsilon$, the absolute value of the integrand in \eqref{eq.main3}, extended by zero to $(0,\infty)\times W$, is bounded by
\begin{equation}\label{eq.dom}
 C\epsilon^{-q}\|\kappa_1\|_\infty\|\kappa_2\|_\infty
 e^{-(\delta-\eta)s}
 \mathbf1_{\{(s,w)\in(0,\infty)\times W:\|w\|^2<3s\}}(s,w),
\end{equation}
where $C$ is independent of $T$ and $\epsilon$.
 Lemma \ref{lem.ratio} gives
 $M_T(s)\leq e^{\eta(s+1)}$ and $M_T(s)\to1$ as $T\to\infty$ for each fixed $s$.
Since $\|w\|^2<3s$ for all $(s,w)\in R_T$, the claim follows from \eqref{eq.sm3}.
 
Since $u\in\operatorname{int}\mathcal C_0$, the cone condition in $R_T$ holds for each fixed $(s,w)$ and all sufficiently large $T$. 
Thus, for every $(s,w)\in(0,\infty)\times W$ outside $\{(s,w)\in(0,\infty)\times W:\|w\|^2=2s\}$, $\mathbf1_{R_T}(s,w)\to\mathbf1_{\{\|w\|^2<2s\}}$.
The exceptional set has product measure zero by Fubini's theorem.
The function in \eqref{eq.dom} is integrable, since its integral in $w$ is a constant multiple of $e^{-(\delta-\eta)s}s^r$. 
Thus \eqref{eq.sm1} and the dominated convergence theorem give
\begin{align*}
 \lim_{T\to\infty}I'_{T,\epsilon}
 &=\nu(\kappa_1)\nu_{\mathrm i}(\kappa_2)
 \int_WJ_u(w)\int_{\|w\|^2/2}^\infty e^{-\delta s}\,ds\,dw
 =c_u\nu(\kappa_1)\nu_{\mathrm i}(\kappa_2).
\end{align*}
Together with \eqref{eq.main2}, this proves (1).
 
For (2), Lemma \ref{lem.ratio}(3) gives, for $0\leq s\leq T/2$,
\begin{equation*}
 |M_T(s)-1|\ll(s+1)e^{\eta(s+1)}T^{-d}+s/T.
\end{equation*}
By \eqref{eq.sm4}, for $(s,w)\in R_T$, we have
\begin{gather*}
 \Psi(T-s)\mathsf S_\epsilon((T-s)u+\sqrt{T-s}\,w)
 =J_u(w)\nu(\kappa_1)\nu_{\mathrm i}(\kappa_2)+E_T(s,w),
\end{gather*}
where
$$
 |E_T(s,w)|
 \ll\left[\epsilon^{-2(q+k)}T^{-d}
 +\epsilon^{-q-1}e^{-cT}\right]
 \|\kappa_1\|_{C^k}\|\kappa_2\|_{C^k}.
$$
For each fixed $s$, $\{w\in W:(s,w)\in R_T\} \subset\{w\in W:\|w\|^2<3s\},$ and its volume is $O(s^r)$.
Substituting these bounds into \eqref{eq.main3} and using boundedness of $J_u$ gives
\begin{multline}\label{eq.main4}
 I'_{T,\epsilon}
 =\nu(\kappa_1)\nu_{\mathrm i}(\kappa_2)
 \int_0^{T/2}\int_W\mathbf1_{R_T}(s,w)e^{-\delta s}J_u(w)\,dw\,ds\\*
 +O\!\left(\left[T^{-1}+\epsilon^{-2(q+k)}T^{-d}
 +\epsilon^{-q-1}e^{-cT}\right]
 \|\kappa_1\|_{C^k}\|\kappa_2\|_{C^k}\right).
\end{multline}
 It remains to bound the
 difference between this integral and $c_u$.
Set
$$
 R_\infty=\{(s,w)\in(0,\infty)\times W:\|w\|^2<2s\}.
$$
Then for each fixed $0<s\leq T/2$,
$$
 \{w\in W:(s,w)\in R_T\setminus R_\infty\} \subset\{w\in W:2s\leq\|w\|^2<2s+2s^2/T\}.
$$
For $q>1$, its volume is at most
$$
 \omega_{q-1}(2s)^r\bigl((1+s/T)^r-1\bigr) \ll T^{-1}s^{r+1},
$$
since $0<s/T\leq1/2$.
When $q=1$, this set is empty.
 
Choose $a>0$ such that $B(u,a)\subset\mathcal C_0$.
If $(s,w)\in R_\infty\setminus R_T$ and $s\leq T/2$, then $(T-s)u+\sqrt{T-s}\,w\notin\mathcal C_0$, so $a^2(T-s)\leq\|w\|^2<2s$.
Hence $s>c_0T$, where $c_0:=a^2/(2+a^2)>0$.
Together with the boundedness of $J_u$, it follows that
\begin{multline*}
 \left|\int_0^{T/2}\int_W
 \mathbf1_{R_T}(s,w)e^{-\delta s}J_u(w)\,dw\,ds-c_u\right|\\* \ll T^{-1}\int_0^\infty e^{-\delta s}s^{r+1}\,ds +\int_{c_0T}^\infty e^{-\delta s}s^r\,ds +\int_{T/2}^\infty e^{-\delta s}s^r\,ds \ll T^{-1}.
\end{multline*}
Combine this with \eqref{eq.main4}.
The exponentially small errors in \eqref{eq.main1} and \eqref{eq.main2} are bounded by $\epsilon^{-2(q+k)}T^{-d}$ for all large $T$, independently of $\epsilon$, since $0<\epsilon\leq\epsilon_0\leq1$ and $2(q+k)\geq q$. 
This proves (2).
\end{proof}

\begin{theorem}\label{thm.eqd}
For all $\kappa_1\in C(\mathcal F_\theta)$ and $\kappa_2\in C(\mathcal F_{\mathrm i\theta})$,
\begin{equation}\label{eq.e12}
 \lim_{T\to\infty}\frac{\Psi(T)e^{-\delta T}}{T^r}
 \sum_{\substack{\gamma\in\Gamma, \|\mu_\theta(\gamma)\|<T}}
 \kappa_1(x_\gamma)\kappa_2(y_\gamma)
 =c_u\nu(\kappa_1)\nu_{\mathrm i}(\kappa_2).
\end{equation}
Hence, for Borel sets $E_1\subset\mathcal F_\theta$ and $E_2\subset\mathcal F_{\mathrm i\theta}$ with $\nu(\partial E_1)=\nu_{\mathrm i}(\partial E_2)=0$,
$$
 \lim_{T\to\infty}\frac{\Psi(T)e^{-\delta T}}{T^r}
 \#\{\gamma\in\Gamma:\|\mu_\theta(\gamma)\|<T,\
 x_\gamma\in E_1,\ y_\gamma\in E_2\}
 =c_u\nu(E_1)\nu_{\mathrm i}(E_2).
$$
\end{theorem}
 
\begin{proof}
Suppose first that $\kappa_1,\kappa_2$ are smooth and nonnegative, and let $L=c_u\nu(\kappa_1)\nu_{\mathrm i}(\kappa_2)$.
By Lemma \ref{lem.ratio}(1), $\lambda_T/\lambda_{T+a}\to e^{\delta a}$ for every fixed $a\in\mathbb R$.
Hence Lemmas \ref{lem.ntit} and \ref{lem.main}(1) give
$$
 e^{-4\delta\epsilon}L
 \leq\liminf_{T\to\infty}\lambda_T\mathcal N_{\kappa_1,\kappa_2}(T)
 \leq\limsup_{T\to\infty}\lambda_T\mathcal N_{\kappa_1,\kappa_2}(T)
 \leq e^{4\delta\epsilon}L.
$$
Letting $\epsilon\to0$ proves \eqref{eq.e12} for  nonnegative smooth $\kappa_1,\kappa_2$.
Taking $\kappa_1=\kappa_2=1$ proves Theorem \ref{thm.cnt1}(1) and bounds the normalized total count.
Linearity and uniform approximation  extend \eqref{eq.e12} to all continuous $\kappa_1,\kappa_2$. 
The claim for Borel sets follows by a standard approximation argument.
\end{proof}
 
 \begin{proof}[Proof of Theorem \ref{thm.cnt2}]
We use the assumed linear $\theta$-regularity near $u$ when applying Lemma \ref{lem.main}(2) and the error term of Theorem \ref{thm.tubes}.
For $|a|\leq2\epsilon_0$, Lemma \ref{lem.ratio}(3) gives
\begin{equation}\label{eq.r1}
 \frac{\lambda_T}{\lambda_{T+a}}
 =e^{\delta a}\left(1+O(T^{-d}+T^{-1})\right).
\end{equation}
Suppose first that $\kappa_1,\kappa_2$ are smooth and nonnegative, and let $L=c_u\nu(\kappa_1)\nu_{\mathrm i}(\kappa_2)$.
Applying Lemma \ref{lem.main}(2) at $T\pm2\epsilon$ and using \eqref{eq.r1} and \eqref{eq.ntit}, we obtain
\begin{multline*}
 |\lambda_T\mathcal N_{\kappa_1,\kappa_2}(T)-L|\\*
 \ll\left[\epsilon+T^{-1}+\epsilon^{-2(q+k)}T^{-d}
 +\epsilon^{-q-1}e^{-cT}\right]
 \|\kappa_1\|_{C^k}\|\kappa_2\|_{C^k},
\end{multline*}
uniformly for $0<\epsilon\leq\epsilon_0$ and all large $T$.
 Writing $\kappa_i$ as in the
 proof of Theorem \ref{thm.tubes} gives the same bound for arbitrary smooth
 $\kappa_1,\kappa_2$.
Choosing $\epsilon=T^{-d/(2(q+k)+1)}$ proves the claim for balls with $\eta=\min\{1,d/(2(q+k)+1)\}$.
The two local limits, with their polynomial errors, were proved in Theorem \ref{thm.llt}.
 
For the claim on the tube count, write $R$ for the radius $\epsilon$ in Theorem \ref{thm.cnt1}(2),
and fix $0<h<R/2$.
Since $t_\gamma\geq0$ and $\|tu+z\|^2=t^2+\|z\|^2$, the count in $\mathbb T_R\cap B(0,T)$ lies, for $T>R$, between the counts in
$$
 [0,T-R^2/T]u+\{z\in W:\|z\|\leq R-h\}
 \quad\text{and}\quad
 [0,T]u+\{z\in W:\|z\|\leq R\}.
$$
Here $(T-R^2/T)^2+R^2<T^2$.
The uniform  bound at the end of the proof of Theorem \ref{thm.tubes} applies to both tubes, with $b_T=-R^2/T$ for the first.
Their normalized main terms differ by $O_R(h+T^{-1})$, and their errors are $O_R(T^{-d/(2(q+k)+1)})$, with the same $k$ as above.
Taking $h=T^{-d/(2(q+k)+1)}$ proves the tube claim with the same $\eta$.
When $W=\{0\}$, the same argument applies, with no need for $h$.
\end{proof}
 
 \section{Counting for the sup norm}\label{sec.sup}
Replacing the Euclidean norm by another norm can change the polynomial factor in the counting asymptotic. 
We illustrate this with norms of the form $\max_j|\ell_j|$, such as the sup norm.
 
Set $q=\dim\mathfrak a_\theta$ and fix linear forms $\ell_1,\ldots,\ell_m$ spanning $\mathfrak a_\theta^*$. 
Set
$$
 N_\infty(v):=\max_{1\leq j\leq m}|\ell_j(v)|,
 \qquad
 \delta_\infty:=\max_{\substack{v\in\mathfrak a_\theta^+, 
                              N_\infty(v)=1}}
                         \psi_\Gamma^\theta(v)>0.
$$
 
Choose $u\in\mathfrak a_\theta^+$ with $N_\infty(u)=1$ and $\psi_\Gamma^\theta(u)=\delta_\infty$.
By \cite[Lemma 3.3.4]{Quint1}, there exists a tangent form $\psi\in\mathfrak a_\theta^*$ at $u$ satisfying
\begin{equation}\label{eq.cond}
 \psi_\Gamma^\theta\leq\psi\leq\delta_\infty N_\infty,
 \qquad \psi(u)=\delta_\infty.
\end{equation}
 
We assume that $\psi>0$ on $\mathcal L_\theta\setminus\{0\}$ and that the supporting hyperplane $\psi=\delta_\infty$ meets the unit ball of $N_\infty$ only at $u$:
\begin{equation}\label{eq.sing}
 \{v\in\mathfrak a_\theta:N_\infty(v)\leq1,
                              \ \psi(v)=\delta_\infty\}
 =\{u\}.
\end{equation}
The positivity assumption implies that $\psi$ is $(\Gamma,\theta)$-proper.
Let $\nu=\nu_\psi$ and $\nu_{\mathrm i}=\nu_{\psi\circ\mathrm i}$ be PS probability measures, and assume that their BMS measure satisfies Definition \ref{lmp} at $u$, with functions $\Psi,J_u$.
Lemma \ref{lem.lm-interior} gives $u\in\operatorname{int}\mathcal L_\theta\subset\operatorname{int}\mathfrak a_\theta^+$.
We apply the results of Section \ref{sec.pc} with these choices and set $W=\ker\psi$ and $\mathfrak j_u=\delta_\infty/\|\psi\|_*$.
 
\begin{theorem}\label{thm.sup}
Under these assumptions, there exists  $c_{\infty}>0$ such that, as $T\to\infty$,
\begin{equation}\label{eq.sup2}
 \#\{\gamma\in\Gamma:N_\infty(\mu_\theta(\gamma))<T\}
 \sim c_{\infty}
       \frac{e^{\delta_\infty T}}{\Psi(T)}.
\end{equation}
 
\end{theorem}
 
The factor $T^{(q-1)/2}$ in the Euclidean ball asymptotic is absent here. 
In particular, when $\Psi(T)=T^{(q-1)/2}$, the sup-norm count has the main term $c_{\infty}T^{-(q-1)/2}e^{\delta_\infty T}$.
The  different powers of $T$ already occur in the integral $\int_{N(v)<T}e^{\psi(v)}\,dv$, where $N$  is the chosen norm.
For the Euclidean norm and the form $\psi$ in Section \ref{sec.eucl},  this integral has order $T^{(q-1)/2}e^{\delta T}$.
For the sup norm under \eqref{eq.sing}, it has order $e^{\delta_\infty T}$.
To see the difference, write $v=tu+z$ with $z\in\ker\psi$ and keep $T-t$ bounded.
The Euclidean ball allows $\|z\|$ of order $\sqrt T$,  giving a factor $T^{(q-1)/2}$  after integration in $z$.
For the sup norm under \eqref{eq.sing}, $z$ remains bounded.
Local mixing  gives the factor $1/\Psi(T)$ and the  function $J_u$ in the leading constant.
The following figure illustrates this in dimension two.
 
\begin{figure}[H]
\centering
\definecolor{sectiongreen}{RGB}{22,128,87}
\definecolor{boundarygray}{RGB}{65,72,72}
\definecolor{ballfill}{RGB}{247,249,248}
\resizebox{0.705\linewidth}{!}{%
\begin{tikzpicture}[x=0.66cm,y=0.66cm,
  every node/.style={font=\small},
  axis/.style={->,line width=0.5pt,black!65},
  boundary/.style={line width=0.85pt,boundarygray},
  ray/.style={->,dashed,line width=0.55pt,black!55},
  cut/.style={dashed,line width=0.5pt,sectiongreen!55},
  chord/.style={sectiongreen,line width=2.0pt,line cap=round}]
\pgfmathsetmacro{\T}{6}
\pgfmathsetmacro{\discB}{\T/sqrt(2)}
\pgfmathsetmacro{\discC}{(\T-1)/sqrt(2)}
\pgfmathsetmacro{\discLo}{(\T-1-sqrt(2*\T-1))/sqrt(2)}
\pgfmathsetmacro{\discHi}{(\T-1+sqrt(2*\T-1))/sqrt(2)}
\pgfmathsetmacro{\squareC}{\T-1/sqrt(2)}
\pgfmathsetmacro{\squareLo}{\T-sqrt(2)}
 
\begin{scope}
  \node[font=\normalsize] at (3.1,7.35) {Euclidean norm};
  \fill[ballfill] (0,0) -- (\T,0)
      arc[start angle=0,end angle=90,radius=6] -- cycle;
  \draw[boundary] (\T,0)
      arc[start angle=0,end angle=90,radius=6];
  \draw[axis] (0,0) -- (7.05,0) node[right] {$x$};
  \draw[axis] (0,0) -- (0,7.05) node[above] {$y$};
  \draw[ray] (0,0) -- (6.62,6.62)
      node[above left,inner sep=2pt] {$x=y$};
  \node[below left,inner sep=3pt] at (0,0) {$0$};
  \node[below,inner sep=3pt] at (\T,0) {$T$};
  \node[left,inner sep=3pt] at (0,\T) {$T$};
  \draw[cut] ({\discLo-0.40},{\discHi+0.40}) --
             ({\discHi+0.40},{\discLo-0.40});
  \draw[chord] (\discLo,\discHi) -- (\discHi,\discLo);
  \fill[sectiongreen] (\discLo,\discHi) circle[radius=1.6pt];
  \fill[sectiongreen] (\discHi,\discLo) circle[radius=1.6pt];
  \draw[<->,line width=0.6pt] (\discC,\discC) -- (\discB,\discB)
     node[midway,inner sep=1.2pt,fill=ballfill] {$1$};
  \fill[boundarygray] (\discB,\discB) circle[radius=1.5pt];
  \node[text=sectiongreen,font=\normalsize] at (3.1,-1.0)
       {$\text{Length}=2\sqrt{2T-1}\asymp\sqrt T$};
\end{scope}
 
\begin{scope}[xshift=6.25cm]
  \node[font=\normalsize] at (3.1,7.35) {Sup norm};
  \fill[ballfill] (0,0) rectangle (\T,\T);
  \draw[boundary] (0,\T) -- (\T,\T) -- (\T,0);
  \draw[axis] (0,0) -- (7.05,0) node[right] {$x$};
  \draw[axis] (0,0) -- (0,7.05) node[above] {$y$};
  \draw[ray] (0,0) -- (6.62,6.62)
      node[above left,inner sep=2pt] {$x=y$};
  \node[below left,inner sep=3pt] at (0,0) {$0$};
  \node[below,inner sep=3pt] at (\T,0) {$T$};
  \node[left,inner sep=3pt] at (0,\T) {$T$};
  \draw[cut] ({\squareLo-0.40},{\T+0.40}) --
             ({\T+0.40},{\squareLo-0.40});
  \draw[chord] (\squareLo,\T) -- (\T,\squareLo);
  \fill[sectiongreen] (\squareLo,\T) circle[radius=1.6pt];
  \fill[sectiongreen] (\T,\squareLo) circle[radius=1.6pt];
  \draw[<->,line width=0.6pt] (\squareC,\squareC) -- (\T,\T)
     node[midway,inner sep=1.2pt,fill=ballfill] {$1$};
  \fill[boundarygray] (\T,\T) circle[radius=1.5pt];
  \node[text=sectiongreen,font=\normalsize] at (3.1,-1.0)
       {$\text{Length}=2$};
\end{scope}
\end{tikzpicture}
} \caption{Sections of the Euclidean and sup-norm balls.}
\label{fig.norm}
\end{figure}
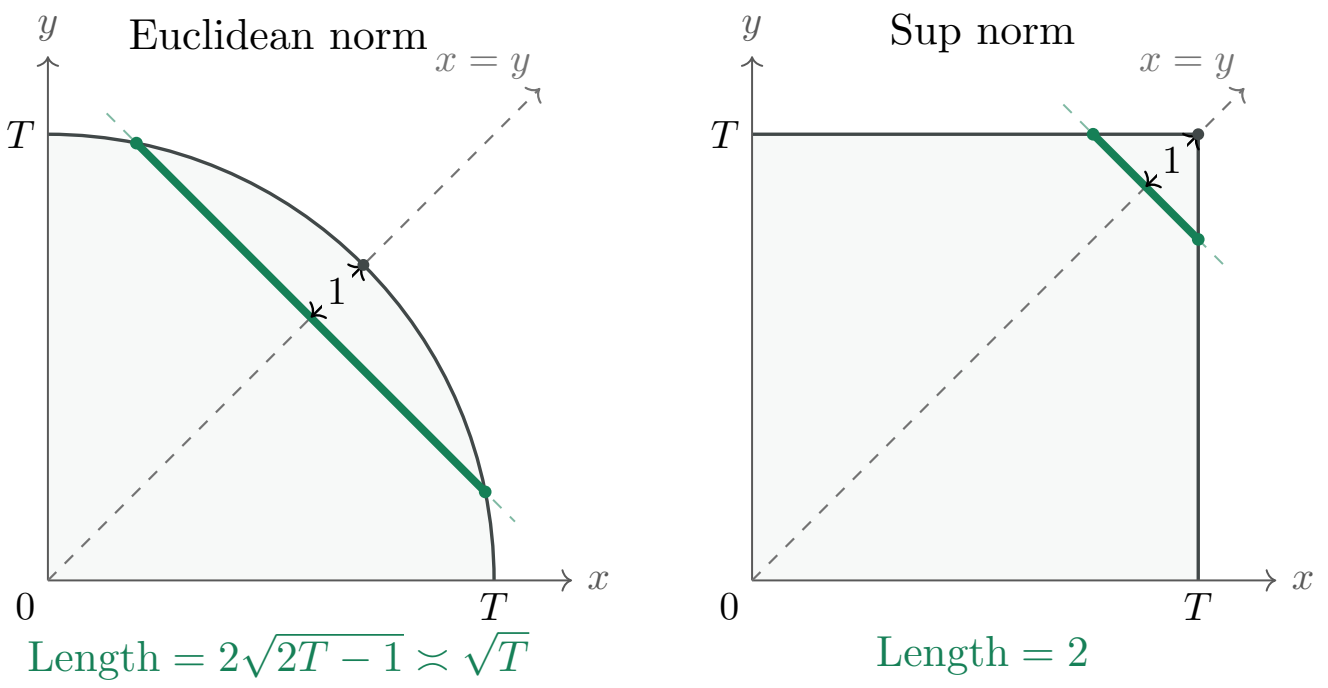

\begin{proof}
By \eqref{eq.cond} and \eqref{eq.sing}, $u$ is the unique maximizer of $\psi_\Gamma^\theta$ on $\{v\in\mathfrak a_\theta^+:N_\infty(v)=1\}$.
Define for $z\in W$,
$$
 p(z):=\max_{\{j:\,|\ell_j(u)|=1\}}\ell_j(u)\ell_j(z).
$$
Let $I:=\{j:|\ell_j(u)|=1\}$.
Since $|\ell_j(u)|<1$ for $j\not\in I$, in a sufficiently small neighborhood of $u$ the maximum defining $N_\infty$ is attained among the indices in $I$. 
Moreover, for $j\in I$, the sign of $\ell_j$ in this neighborhood agrees with that of $\ell_j(u)$.
After shrinking $\mathcal C_0$ around $u$ if necessary, these properties hold throughout $\mathcal C_0\setminus\{0\}$.
We also arrange that $\psi>0$ on $\mathcal C_0\setminus\{0\}$, since $\psi(u)>0$.
Since $\ell_j(u)^2=1$ for $j\in I$, we obtain
$$
 N_\infty(tu+z) =\max_{j\in I}\bigl(t+\ell_j(u)\ell_j(z)\bigr) =t+p(z)
$$
for all $t>0$ and $z\in W$ satisfying $tu+z\in\mathcal C_0$.
Applying $\psi\leq\delta_\infty N_\infty$ to $u+az$ for sufficiently small $a>0$ gives $p(z)\geq0$ for every $z\in W$.
If $p(z)=0$ for some nonzero $z\in W$, then, for sufficiently small $a>0$, $N_\infty(u+az)=1$ and $\psi(u+az)=\delta_\infty$, contrary to \eqref{eq.sing}.
By homogeneity and compactness, there is $c>0$ such that $p(z)\geq c\|z\|$ for all $z\in W$.

Take $\kappa_1=\kappa_2=1$ in $\mathsf S_\epsilon$ and set
$$
 I^\infty_{T,\epsilon} :=\int_{N_\infty(v)<T}e^{\psi(v)} \mathsf S_\epsilon(v)\,dv.
$$
Since $u$ is the unique maximizing vector, the proofs of Lemmas \ref{lem.subexp} and \ref{lem.subexp2} apply with $N_\infty$ in place of the Euclidean norm, using
$\psi_\Gamma^\theta\leq\psi$ and \eqref{eq.sing} in place of Lemma \ref{lem.dir}. 
Thus, for fixed $\epsilon$, restricting this integral to $v=tu+z\in\mathcal C_0$ with $t\geq t_0$ changes its product with $\Psi(T)e^{-\delta_\infty T}$ by $o(1)$.
 
Fix $0<\eta<\delta_\infty$ and choose $t_0$ sufficiently large for \eqref{eq.sm3} and Lemma \ref{lem.ratio}(2).
Since $\psi>0$ on $\mathcal C_0\setminus\{0\}$, the section $\{z\in W:tu+z\in\mathcal C_0\}$ is bounded and has volume $O(t^{q-1})$. 
These  bounds show that the normalized contribution from $t_0\leq t<T/2$ is
$$
 O\!\left( \epsilon^{-q}T^q e^{-(\delta_\infty-\eta)T/2} \right)=o(1).
$$
In the remaining integral, set $s=T-t$ and use $dv=\mathfrak j_u\,dt\,dz$. 
Let
$$
 D_T:=\{(s,z)\in(0,T/2]\times W: (T-s)u+z\in\mathcal C_0,\ p(z)<s\}.
$$
Then
\begin{align*}
 \Psi(T)e^{-\delta_\infty T}I^\infty_{T,\epsilon}
 &=\mathfrak j_u
 \int_0^{T/2}\int_W
 \mathbf1_{D_T}(s,z)e^{-\delta_\infty s}
 \frac{\Psi(T)}{\Psi(T-s)}\\
 &\qquad\cdot
 [\Psi(T-s)\mathsf S_\epsilon((T-s)u+z)]
 \,dz\,ds+o(1).
\end{align*}
For each fixed $(s,z)$, the indicator tends to $\mathbf1_{\{p(z)<s\}}$.
By \eqref{eq.sm2} and Lemma \ref{lem.ratio}(1), the remaining factors tend to $e^{-\delta_\infty s}J_u(0)$.
Moreover, \eqref{eq.sm3}, applied with $w=z/\sqrt{T-s}$, and Lemma \ref{lem.ratio}(2) bound the integrand by a constant multiple of
$$
 \epsilon^{-q}e^{-(\delta_\infty-\eta)s} \mathbf1_{\{c\|z\|<s\}}.
$$
This is integrable on $(0,\infty)\times W$.
Dominated convergence therefore gives
\begin{align*}
 \lim_{T\to\infty}
 \Psi(T)e^{-\delta_\infty T}I^\infty_{T,\epsilon}
 &=\mathfrak j_uJ_u(0)
   \int_W\int_{p(z)}^\infty e^{-\delta_\infty s}\,ds\,dz\\
 &=\frac{\mathfrak j_uJ_u(0)}{\delta_\infty}
   \int_W e^{-\delta_\infty p(z)}\,dz
 =:c_{\infty}.
\end{align*}
The bound $p(z)\geq c\|z\|$ shows that this constant is finite, and it is positive.
 
Finally, \eqref{eq.ite} and the triangle inequality for $N_\infty$ give constants $a,b>0$ such that
$$
 e^{-b\epsilon}I^\infty_{T-a\epsilon,\epsilon} \leq \#\{\gamma\in\Gamma:N_\infty(\mu_\theta(\gamma))<T\} \leq e^{b\epsilon}I^\infty_{T+a\epsilon,\epsilon}.
$$
Multiply by $\Psi(T)e^{-\delta_\infty T}$ and apply Lemma \ref{lem.ratio}(1). 
Letting $T\to\infty$ and then $\epsilon\to0$ proves \eqref{eq.sup2}.
\end{proof}
\subsection*{The sup norm in $\operatorname{SL}_3(\mathbb R)$}
Let $G=\operatorname{SL}_3(\mathbb R)$ and $\theta=\Pi$, and let
\begin{gather*}
 N_\infty(\operatorname{diag}(v_1,v_2,v_3)):=\max_i|v_i|,\\
 u_\infty:=\operatorname{diag}(1,0,-1),\qquad
 Z:=\operatorname{diag}(1,-2,1).
\end{gather*}
Opposition is $\mathrm i(v_1,v_2,v_3)=(-v_3,-v_2,-v_1)$, so
$$
 \frac{v+\mathrm i v}{2}=\frac{v_1-v_3}{2}u_\infty.
$$
Choose $v\in\operatorname{int}\mathcal L_\Pi$. Since $\mathcal L_\Pi$ is convex and $\mathrm i$-invariant,
this midpoint lies in $\operatorname{int}\mathcal L_\Pi$, and hence $u_\infty\in\operatorname{int}\mathcal L_\Pi$.
Set $\delta_\infty:=\psi_\Gamma^\Pi(u_\infty)>0$.  For every $v\in\mathcal L_\Pi$, concavity and opposition invariance give
$$\psi_\Gamma^\Pi(v)
 \leq\psi_\Gamma^\Pi\left(\frac{v+\mathrm i v}{2}\right)
 =\delta_\infty\frac{v_1-v_3}{2}
 \leq\delta_\infty N_\infty(v).
$$
Thus $\psi_\infty(v):=\delta_\infty(v_1-v_3)/2$ is tangent at $u_\infty$.
 If $N_\infty(v)\leq1$ and $\psi_\infty(v)=\delta_\infty$,
 then $v_1=1$ and $v_3=-1$, and the trace-zero condition gives $v_2=0$. This proves \eqref{eq.sing}
 and the uniqueness of the maximizing vector.
The form $\psi_\infty$ is positive on $\mathfrak a^+\setminus\{0\}$,  so it is proper.
Moreover, $W=\mathbb RZ$, and since $I=\{1,3\}$,
$p(aZ)=\max\{a,-a\}=|a|$.
Whenever the associated BMS measure satisfies the local mixing hypothesis at $u_\infty$,
Theorem \ref{thm.sup} therefore applies.
 
\begin{proof}[Proof of Corollary \ref{cor.sup}]
Since $u_\infty\in\operatorname{int}\mathcal L_\Pi$, the required local mixing with $\Psi(t)=\sqrt t$ follows from \cite[Theorems 1.1(2) and 17.1(2)]{KOP}.
Applying Theorem \ref{thm.sup} proves \eqref{eq.sup}.
\end{proof}

\end{document}